\documentclass[12pt]{article}

\usepackage{amsmath}
\usepackage{amssymb}
\usepackage[curve,matrix,arrow,cmtip]{xy}
\NoComputerModernTips

\newtheorem{theorem}{\indent Theorem}[section]
\newtheorem{lemma}[theorem]{\indent Lemma}
\newtheorem{prop}[theorem]{\indent Proposition}
\newtheorem{cor}[theorem]{\indent Corollary}

\newtheorem{remark}[theorem]{\indent Remark}

\newcommand{\rhobar}{{\overline{\rho}}}

  {\vskip 0.3ex plus 0.5ex minus 0ex \pagebreak[1]
   \trivlist
   \item[\hskip\labelsep \textbf{General Convention}]}%
  {\endtrivlist
   \vskip .3ex plus 1ex minus 0ex \pagebreak[2]}

\advance\textheight by .06cm

\long\def\forget#1{}

\def\?{\ ???\ \immediate\write16{}%
\immediate\write16{Warning: There was still a question mark . . . }%
\immediate\write16{}}

\def\Trace{\mathrm{Trace}}
\DeclareMathOperator{\image}{im}

\def\unr{{\mathrm{unr}}}
\def\BF{\mathbf{F}}

\def\BQ{\mathbf{Q}}
\def\BN{\mathbf{N}}

\def\BZ{\mathbf{Z}}
\def\eps{\varepsilon}
\def\ad{\mathrm{ad}^0(\rhobar)}

\def\Ad{\mathrm{ad}(\rhobar)}

\def\GL{\mathrm{GL}}

\def\End{\mathrm{End}}

\def\SL{\mathrm{SL}}

\def\Frob{\mathrm{Frob}}
\def\Gal{\mathrm{Gal}}

\def\Hom{\mathrm{Hom}}

\def\rank{\mathrm{rank}}

\DeclareMathOperator{\mult}{mult}

\def\CA{\mathcal{A}}
\def\CM{\mathcal{M}}

\def\CC{\mathcal{C}}

\def\CK{\mathcal{K}}
\def\CS{\mathcal{S}}
\def\CT{\mathcal{T}}

\def\Fm{\mathfrak{m}}

\def\theenumi{(\roman{enumi})}

\def\p@enumii{\theenumi}

\def\longto{\longrightarrow}
\def\into{\hookrightarrow}

\newbox\mybox

\def\arrover#1{\mathrel{
       \setbox\mybox=\hbox spread 1.4em{\hfil$\scriptstyle#1$\hfil}
       \vbox{\offinterlineskip\copy\mybox
             \hbox to\wd\mybox{\rightarrowfill}}}}
\def\larrover#1{\mathrel{
       \setbox\mybox=\hbox spread 1.4em{\hfil$\scriptstyle#1$\hfil}
       \vbox{\offinterlineskip\copy\mybox
             \hbox to\wd\mybox{\leftarrowfill}}}}

\def\ontoover#1{\mathrel{
       \setbox\mybox=\hbox spread 1.4em{\hfil$\scriptstyle#1$\hfil}
       \vbox{\offinterlineskip\copy\mybox
             \hbox to\wd\mybox{\rightarrowfill\hskip-2.8mm
                               $\rightarrow$}}}}
\def\leftontoover#1{\mathrel{
       \setbox\mybox=\hbox spread 1.4em{\hfil$\scriptstyle#1$\hfil}
       \vbox{\offinterlineskip\copy\mybox
             \hbox to\wd\mybox{$\leftarrow$\hskip-2.8mm
                               \leftarrowfill}}}}
\def\onto{\ontoover{\ }}

\newif\ifnormalesBeweisEnde

\newenvironment{Proofof}[1]%
  {\vskip 0.3ex plus 0.5ex minus 0ex \pagebreak[1]
   \global\normalesBeweisEndetrue
   \trivlist
   \item[\hskip\labelsep {\textsc{Proof} \rm of #1}:]}%
  {\ifnormalesBeweisEnde \EndOfBeweis \fi
   \endtrivlist
   \vskip 1ex plus 1ex minus 0ex \pagebreak[2]}
\newenvironment{Proof}%
  {\vskip 0.3ex plus 0.5ex minus 0ex \pagebreak[1]
   \global\normalesBeweisEndetrue
   \trivlist
   \item[\hskip\labelsep \textsc{Proof}:]}%
  {\ifnormalesBeweisEnde \EndOfBeweis \fi
   \endtrivlist
   \vskip 1ex plus 1ex minus 0ex \pagebreak[2]}
\def\EndOfBeweis{\hskip .5em \vrule width 1.0ex height 1.0ex depth 0.3ex}

\newcommand{\notdiv}{\mathopen{\mathchoice
             {\not{|}\,}
             {\not{|}\,}
             {\!\not{\:|}}
             {\not{|}}
             }}

\begin{document}
%

\hfuzz=2pt

\title{Mod $\ell$ representations of arithmetic fundamental groups I \\ 
{\small (An analog of Serre's conjecture for function fields)}}

\author{Gebhard B\"ockle \thanks{G.B.\ would like to thank the
TIFR for its hospitality in the Summer of 2002 during which the first
decisive steps in this work were made and also the ETH for its
inspirational environment and for its generosity in providing a three
year post-doctoral position} \\ 
Chandrashekhar Khare \thanks{{some of the work on this paper was done during a visit to Universit\'e Paris 7 which was supported 
by Centre franco-indien pour la promotion de la recherche avanc{\'e}e (CEFIPRA) under
Project 2501-1 {\it Algebraic Groups in Arithmetic and Geometry}\newline
\indent 2000 MSC: 11F80, 11F70, 14H30, 11R34}}}

\date{}

\maketitle

\begin{abstract}
There is a well-known conjecture of Serre that any continuous, irreducible representation
$\rho:G_{\bf Q} \rightarrow GL_2(\overline{{\bf F}}_p)$ which is odd arises from a newform.
Here $G_{\bf Q}$ is the absolute Galois group of ${\bf Q}$ and $\overline{{\bf F}}_p$ an algebraic closure 
of the  finite field ${\bf F}_p$ of $p$ elements. We formulate such a conjecture for $n$-dimensional mod $\ell$ representations
of $\pi_1(X)$, for any positive integer $n$, and where $X$ is a geometrically irreducible, smooth curve over a finite
field $k$ of characteristic $p$ ($p \neq \ell$), and prove this conjecture in a large number of cases.  In fact a proof of all cases
of the conjecture for 
$\ell>2$ follows from a result announced 
(conditionally) by Gaitsgory in \cite{gaitsgory}.  The methods are completely different.
\end{abstract}

\section{Statement of main result: analog of Serre's conjecture for function fields}
Let $X$ be a geometrically irreducible, smooth curve over a finite
field $k$ of characteristic $p$ and cardinality $q$. Denote by 
$K$ its function field and by $\widetilde{X}$ its smooth
compactification and set $S:=\widetilde X\setminus X$. 
Let $\pi_1(X)$ denote the arithmetic fundamental group of $X$. 
Thus $\pi_1(X)$ sits in the exact sequence $$0 \rightarrow
\pi_1(\overline{X}) \rightarrow \pi_1(X) \rightarrow G_k \rightarrow 0,$$
where $\overline{X}$ is the base change of $X$ to an algebraic closure
of $k$, and $G_F$ denotes the absolute Galois group of any field~$F$.

We study here mod $\ell$ representations of $\pi_1(X)$, i..e,
continuous absolutely irreducible representations 
$\rhobar:\pi_1(X) \rightarrow \GL_n({\bf F})$ with ${\bf F}$
a finite field of characteristic $\ell \neq p$. In this paper we are mainly
interested in an analog of (the qualitative part of) Serre's conjectures in
\cite{serre} in the function field situation. 

Let us fix once and for all an embedding $\iota\!:\overline{\bf Q}
\hookrightarrow \overline{{\bf Q}}_{\ell}$. 
Then with respect to this embedding $\iota$, there is a correspondence between
$n$-dimensional $\ell$-adic representations of $\pi_1(X\setminus T)$
with finite order determinant and suitably ramified cuspidal
eigenforms (or equivalently cuspidal automorphic representations 
with a newvector fixed by a suitable open compact of $GL_n({\bf A}_K)$) on $GL_n({\bf A}_K)$ with finite central order character,  for
any finite set $T$ of places of $X$. This correspondence 
is the global Langlands correspondence for functions due to
Drinfeld and Lafforgue (see \cite{drinfeld} and \cite{lafforgue}).

We call a residual representation $\rhobar$  {\em automorphic}, if it is isomorphic
to the residual representation  attached to (an
integral model of) an $n$-dimensional continuous representation 
$\pi_1(X\setminus T) \rightarrow GL_n(\overline{ {\bf Q}}_{\ell} )$ 
that is associated to a cuspidal automorphic representation 
of  $GL_n({\bf A}_K)$ in \cite{drinfeld} and \cite{lafforgue} for some finite set of places $T$. An analog of Serre's conjecture in
the function field setting is therefore that any absolutely
irreducible residual representation $\rhobar$ is
automorphic. It is worth noting that unlike in the classical setting here there are no ``local conditions''
that need to be imposed on $\rhobar$ to expect it to be automorphic. In view of \cite{lafforgue} this conjecture is
equivalent to the assertion that any such $\rhobar$ lifts to an $\ell$-adic representation of $\pi_1(X\setminus T)$ of finite order
determinant for some finite subset $T$ of~$X$. 

There is little known about Serre's original conjecture, while the
analog that we study for function fields is more accessible because
of the results in \cite{drinfeld} and \cite{lafforgue}. Moreover the main
result of \cite{dejong} directly implies that the function field analog of
Serre's conjecture holds for $n\leq 2$ (for $n=1$ it is a simple
consequence of class field theory). This is strong evidence in favour
of the analog. In fact, the main conjecture made in \cite{dejong}
may be regarded as a refinement of the above analog and easily
implies~it.

\smallskip

In Theorem~\ref{OnSerreConj} below, we establish the analog in many
more cases by producing suitable $\ell$-adic liftings of
$\rhobar$. Our approach uses the Galois cohomological methods of
R.\ Ramakrishna, \cite{ramakrishna}, and their further refinements by
R.\ Taylor in \cite{taylor}. The following is an important special case
of Theorem~\ref{OnSerreConj}:

\begin{theorem}\label{main}
\label{SpecialSerreConj}
  Let $X$ be a smooth, geometrically irreducible 
  curve defined over a finite field $k$ of
  characteristic $p$, and $\rhobar:\pi_1(X) \rightarrow \SL_n({\bf F})$ be a
  representation with ${\bf F}$ a finite field of characteristic
  $\ell \neq p$. Assume that
  \begin{itemize}
  \item[(i)] $\rhobar$ has full image, $|{\bf F}|\ge4$, $\ell\notdiv
  n$, and
  \item[(ii)] at any $v\in S$ the ramification is
  either tame or of order prime to~$\ell$.
  \end{itemize}
  Then $\rhobar$ lifts to a representation $\rho:\pi_1(X
  \backslash T) \rightarrow \SL_n(W({\bf F}))$ with $T$ a finite set of
  places of $X$ and $W({\bf F})$ the Witt vectors of ${\bf F}$. Hence
  $\rhobar$ is automorphic.
\end{theorem}

What is mainly needed in the proof are
that the adjoint representation $\ad$ of $\rhobar$ on the traceless
matrices of $M_n(\BF)$ is irreducible and that
$H^1(\image(\rhobar),\ad)$ is (almost) zero. It is this lifting
result, which is a consequence of de Jong's conjecture, and that we  
prove under the technical assumption above, that seems
to be crucial for the applications of de Jong's conjecture by Drinfeld
in \cite{drinfeld1} to some purity conjectures of  
Kashiwara on perverse sheaves. For $\ell>2$ a proof of all cases of Serre's 
conjecture, 
that depends on some ``unpublished mathematics'', follows from the work of 
Gaitsgory, cf.\ \cite{gaitsgory}. The methods are completely different and
while Gaitsgory's work should prove the conjecture in totality for $\ell>2$,
our methods apply in characteristic 2.

\medskip

In a continuation (part II) of this work we study the conjecture of
A.J.~de Jong from \cite{dejong} which is about deformations of representations of
the type $\rhobar$ studied in this paper. For this we will use the
lifting result of this paper. In fact proving de Jong's conjecture was
the main motivation for this work. Our results towards  de Jong's
conjecture yields that in many cases $\rhobar$ arises
from a cuspidal eigenform form of level the conductor of $\rhobar$,
where by ``arises from'' we mean is isomorphic to the reduction of the
$n$-dimensional $\ell$-adic representation (which might no longer have
coefficients in Witt vectors) associated to the eigenform, thus proving results towards the analog of Serre's conjecture in its ``quantitative aspect''. 

\section{Proof of the main result}\label{GalM}

Our main goal is to prove a general criterion for
a residual representation to lift to a charactersitic 0 representation which will then give a proof of
Theorem \ref{main}, upon using Lafforgue's theorem. We start by first making
all the necessary definitions to state a result, Theorem 
\ref{OnSerreConj}, that is more general than
Theorem \ref{main}, but is more technical to state. After stating
the Theorem \ref{OnSerreConj}, we will first
quickly derive Theorem~\ref{SpecialSerreConj} from it. 
Then in the following sections we will, following Ramakrishna,
\cite{ramakrishna}, and Taylor, \cite{taylor}, give the proof
of Theorem \ref{OnSerreConj}. For the general background on Galois cohomology of function
fields, the reader is referred to \cite{nsw}, Ch.7 and~8.

\medskip

Let us fix some notation. For a place $v$ of $\widetilde X$ 
denote by $q_v$ the cardinality of the residue field at~$v$.
Let $G_v\supset I_v\supset P_v$ be the absolute Galois group of the
completion of $K$ at $v$, its inertia and wild inertia
subgroup, respectively. We also choose an embedding $G_v\to G_K$.
For any curve $X\subset \widetilde X$ this yields morphisms $G_v/I_v\into
\pi_1(X)$, and by $\Frob_v\in\pi_1(X)$ we denote the corresponding
Frobenius-substitution at~$v\in X$.

We fix a residual representation $\rhobar\!:\pi_1(X)\to\GL_n(\BF)$, where
$\BF$ is some finite field of characteristic $\ell\neq p$. Denote by
$\Ad$ the module $M_n({\bf F})$ considered as a $\pi_1(X)$-module via
the adjoint action composed with $\rhobar$ and by $\ad$ the
subrepresentation on the traceless matrices $M^0_n({\bf F})$ of
$M_n({\bf F})$.
If $M$ is an $\BF[\pi_1(X)]$-module, then $M(i)$, $i\in\BZ$, denotes
the twist of $M$ by the $i$-th tensor power of the cyclotomic mod $\ell$
character $\chi\!:\pi_1(X)\to\BF_\ell^*$, and $M^*$ denotes the
representation $\Hom(M,\BF)$. Note that $\ad\cong\ad^*$ if
$\ell\notdiv n$. For $M$ as above, we also define
$h^i(\pi_1(X),M):=\dim_{\BF}H^i(\pi_1(X),M)$ and similarly with
$\pi_1(X)$ replaced by some~$G_v$.

\medskip

To state the main technical theorem, we need to introduce some
further notation. Let $\zeta_\ell$ be a primitive $\ell$-th root of unity of
$K$ and denote by $E$ the splitting field of $\rhobar$ over $K$, i.e.,
the fixed field of $\rhobar$ in a fixed separable closure of
$K$. Recall that a matrix $A\in\GL_n(\BF)$ is called {\em regular}, if
$\dim_\BF M_n(\BF)^A=n$, where $A$ operates via the adjoint action,
i.e., via conjugation. 

We call a conjugacy class $[\sigma]$ of $\Gal(E(\zeta_\ell)/K)$ an {\em
$R$-class or Ramakrishna-class for $\rhobar$} if $A:=\rhobar(\sigma)$
is regular and one of the following two cases holds: 
\begin{itemize}
\item[(I)] $\chi(\sigma)\neq1$ and $A$ has distinct simple roots
  $\lambda,\lambda'\in\BF$ with $\lambda'=\chi(\sigma)\lambda$. 
\item[(II)] $\chi(\sigma)=1$ and in the Jordan decomposition of $A$
  there occurs at least one $2\times2$ block with eigenvalue
  $\lambda\in\BF$. 
\end{itemize}
We call a place $v$ of $X$ an {\em $R$-place for $\rhobar$} if the
class of $\Frob_v$ in $\Gal(E(\zeta_\ell)/K)$ is an $R$-class. Note
that if an $R$-class exists, then by the \v{C}ebotarev density theorem,
there exist infinitely many $R$-places. 

While the main use of $R$-places is to provide locally some freedom
for lifting, they can also be useful to remove global obstructions. To
describe this, denote by $V$ the space $\BF^n$ considered as a
representation of $\pi_1(X)$ via $\rhobar$ and let $\sigma$ be an
$R$-class of type (II). The indecomposable summands of $V$ are denoted
by $V_i$, and $\Ad_i$ ($\ad_i$) denotes the corresponding representation
on (the trace zero matrices of) $\End(V_i)$, considered as a
representation of $\langle\sigma\rangle$. Let $\lambda_i$ be one of
the eigenvalues of $V_i$. We define
\begin{equation}\label{DefOfAdSig}
\ad_\sigma:=\prod_{{\mult(\lambda_i)=2}\atop{\lambda_i\in\BF}} \ad_i.
\end{equation}
Since $V\cong \oplus_iV_i$ there is a
$\langle\sigma\rangle$-morphism $\ad\to \ad_\sigma$. 
We say that $\rhobar$ {\em admits sufficiently many $R$-classes} if there exists
at least one $R$-class and if the restriction morphisms
(composed with $\ad\to \ad_\sigma$ at $R$-places)
\begin{eqnarray*}
H^1(\Gal(E(\zeta_\ell)/K),\ad)\!\!\!\!\hspace*{-.11pt}&\to&\hspace*{-.11pt}\!\!\!\!\!\!\!\!\!\!\!\!
\prod_{\genfrac{}{}{0pt}{}{\sigma\,\mathrm{an}\,R\mathrm{-class}}{\mathrm{of}\,\mathrm{type(II)}}}\!\!\!\!\!\!\!\!
H^1(\langle\sigma\rangle,\ad_\sigma)
\!\oplus\!\prod_{v\in S}\!H^1(\rhobar(I_v),\ad),\\ 
H^1(\Gal(E(\zeta_\ell)/K),\ad(1))\!\!\!\!\hspace*{-.11pt}&\to&\hspace*{-.11pt}
\!\!\!\!\!\!\!\!\!\!\!\!\prod_{\genfrac{}{}{0pt}{}{\sigma\,\mathrm{an}\,R\mathrm{-class}}{\mathrm{of}\,\mathrm{type(II)}}}\!\!\!\!\!\!\!\!
H^1(\langle\sigma\rangle,\ad_\sigma)
  \!\oplus\!\prod_{v\in S}\!H^1(\rhobar(I_v),\ad(1))\\ 
\end{eqnarray*}
are injective. Note that $\ad^*_\sigma(1)=\ad_\sigma^*=\ad_\sigma$ for
$R$-classes of type~(II).

Our main result in this chapter is
\begin{theorem}\label{OnSerreConj}
  Let $X$ be a smooth, geometrically irreducible 
  curve defined over a finite field $k$ of
  characteristic $p\neq l$, and $\rhobar:\pi_1(X) \rightarrow
  \GL_n({\bf F})$ be a continuous representation. Assume that
  \begin{itemize}
  \item[(a)] $\ad$ is irreducible over $\BF_\ell[\image(\rhobar)]$,
  \item[(b)] $\rhobar$ has sufficiently many $R$-classes,
  \item[(c)] at all $v\in S$ the ramification is
  either tame or of order prime to~$\ell$.
  \end{itemize}
  Then $\rhobar$ lifts to a representation $\rho:\pi_1(X
  \backslash T) \rightarrow \GL_n(W({\bf F}))$ where
\begin{enumerate}
\item $T$ is a finite set of places of $X$,
\item $\det\rho$ is the Teichm\"uller lift of $\det\rhobar$,
\item for $v\in S$ the conductors of $\rho$ and $\rhobar$ agree, and
\item if $\rhobar$ is tame at $v$, then 
   $\rho(I_v)\stackrel\cong\to\rhobar(I_v)$, i.e., $\rho$ is minimal
   at $v$.
\end{enumerate}
\end{theorem}
Note that we do {\em not} need that $\ad$ is absolutely irreducible.
Note also that the condition that $\ad$ is irreducible implies
that $\ell$ does not divide $n$, since in the case $\ell|n$, the
representation $\ad$ contains the trivial representation on scalar
matrices as a non-trivial submodule.

As an application of Lafforgue's theorem, we find.
\begin{cor}
Any $\rhobar$ as in the previous Theorem is automorphic.
\end{cor}
We have the following example for the existence of sufficiently many
$R$-classes. Combined with  the above theorem it completes the proof
of Theorem~\ref{SpecialSerreConj}. 
\begin{prop}\label{SL_nAdmitsRclasses}
Suppose $\rhobar:\pi_1(X) \rightarrow \SL_n({\bf F})$ is surjective,
$\ell\notdiv n$, $\ell\neq p$ and  $|{\bf F}|>4$. Then $\rhobar$ admits
sufficiently many $R$-classes.
\end{prop}
\begin{Proof}
Let us first show the injectivity of the two restriction morphisms
considered above Theorem~\ref{OnSerreConj}.
If $|\BF|>5$ or $n>2$, then by \cite{cps} we have
$H^1(\SL_n(\BF),M^0_n(\BF))=0$. In this case it easily follows, e.g.,
\cite{boeckle2}, \S 5, that 
$$H^1(\Gal(E(\zeta_\ell)/K),\ad)=H^1(\Gal(E(\zeta_\ell)/K),\ad(1))=0.$$

In the remaining case, note first that by loc.\ cit., if $\chi$ is
non-trivial, one has $H^1(\Gal(E(\zeta_\ell)/K),\ad(1))=0$. 
Finally, in \cite{taylor} it is shown for $n=2$ and
$\BF=\BF_5$ how to find an $R$-class $\sigma$ such that the kernel of
$$ \BF\cong H^1(\Gal(E(\zeta_\ell)/K),\ad) \to H^1(\langle\sigma\rangle,\ad)$$
is trivial (in this particular case one has $\ad_\sigma=\ad$). If
$\chi$ is trivial, the same class also works for $\ad(1)=\ad$.

It remains to prove the existence of at least one $R$-class. For this,
note that $\SL_n(\BF)$ has no abelian quotients, and therefore the
morphism $$\rhobar\times\chi\!:\pi_1(X)\longto
\GL_n(\BF)\times\BF$$ surjects onto $\SL_n(\BF)\times\image(\chi)$.
Since $\SL_n(\BF)$ contains matrices of type (II), the existence of an
$R$-class is obvious. Furthermore, if $\image(\chi)$ is non-trivial,
then one may also find matrices of type (I). This completes the proof
of the proposition. 
\end{Proof}

\subsection{Strategy of the proof of Theorem~\ref{OnSerreConj}}\label{Outline}
Our method of producing lifts is essentially that of 
Ramakrishna, cf.\ \cite{ramakrishna}. However we will follow the more
axiomatic treatment as presented in \cite{taylor}. Let us fix from now 
on a representation $\rhobar\!:\pi_1(X)\to\GL_n(\BF)$ which satisfies
the conditions of Theorem~\ref{OnSerreConj} and let~$n\geq2$, since $n=1$ is
trivial by using Teichm\"uller lifts. Therefore in the following we
assume that $\ad$ is irreducible over $\BF_\ell[\image(\rhobar)]$ (and
hence also $\ell\notdiv n$).
Also define $\eta\!:\pi_1(X)\to W(\BF)$ as the Teichm\"uller lift of
$\det\rhobar$ and for any place $v$ define restrictions
$\eta_v:=\eta_{|G_v}$ and $\rhobar_v:=\rhobar_{|G_v}$. 

\smallskip

The strategy in \cite{ramakrishna} to produce lifts of $\rhobar$ 
to $W(\BF)$ is to first consider 
all deformations of $\rhobar$ which are representations of
$\pi_1(X\setminus T)$ for some fixed finite subset $T$ of $R$-places
of $X$ and which at the places in $S\cup T$ are allowed to only have
ramification of a very specific type. Without loss of generality, we
assume that $\rhobar$ is ramified at the places in $S$ and call these
{\em residually ramified places} or simply $r$-places.

The type of ramification is most conveniently
formulated in terms of local deformation problems $\CC_v$ at places
$v\in S\cup T$. In this formulation, locally the crucial 
requirement is that the resulting versal hull is smooth over the ring
$W(\BF)$ of relative dimension $h^0(G_v,\ad)$. In
Sections~\ref{AuxPrimes} 
and~\ref{RamPrimes}, we will define such $\CC_v$ for $R$- and
$r$-places, respectively. 

The global requirement on $T$ and the types $\CC_v$ is made in such a
way that one can inductively construct lifts of $\rhobar$ to the rings
$W_n(\BF)$ of Witt vectors of length $n$. It can be entirely
formulated in terms of Galois cohomology. In this section we will 
recall the necessary background from \cite{taylor} and give a proof of
the main  theorem pending on a key lemma, whose proof will be given in the
later Section~\ref{KeyLemSec}. 

\medskip

Let $\CA$ denote the category of complete noetherian local
$W(\BF)$-algebras $(R,\Fm_R)$ with residue field $\BF$ and where
morphisms are morphisms of local rings which are the identity on the
residue field. By a {\em lift of determinant $\eta_v$} 
of $\rhobar_v$ we mean a continuous 
representation $\rho\!:G_v\to\GL_n(R)$ for some $(R,\Fm_R)\in\CA$ 
such that $\rho\!\!\pmod{\Fm_R}=\rhobar_v$ and $\det\rho=\eta_v$.

We call a pair $(\CC_v,L_v)$, where $\CC_v$ is a collection of
lifts of $\rhobar_v$ of determinant $\eta_v$, and  where $L_v$ is a
subspace of $H^1(G_v,\ad)$, {\em (locally) admissible and compatible
  with $\eta_v$} if it satisfies the conditions
P1--P7 of \cite{taylor}, where in loc.~cit.\ one has to replace $\Fm$
by $\Fm_R$ and $M_2(\Fm)$ by $M_n(\Fm_R)$ in property~P2. 
\begin{remark}\label{Heuristically}
{\em Observe that the conditions of
loc.~cit.\ imply in particular that the versal hull of the
deformation problem $\CC_v$ exists and the corresponding versal
deformation ring is smooth over $W(\BF)$ of relative dimension
$\dim~L_v$. Heuristically one would expect the versal deformation ring
of all deformations of $\rhobar_{v}$ with fixed
determinant to be a complete intersection, flat over $W(\BF)$ and of
relative dimension $h^0(G_v,\ad)$. Therefore one expects
$\dim L_v\leq h^0(G_v,\ad)$.}
\end{remark}

Suppose one is given a finite set  $T\subset X$ and for each $v\in
S\cup T$ a locally admissible pair $(\CC_v,L_v)$ compatible with 
$\eta_v$. Then a {\em lift of type
$(\CC_v)_{v\in S\cup T}$}, 
is a continuous representation $\rho\!:\pi_1(X\setminus T)\to\GL_n(R)$
for some $(R,\Fm_R)\in\CA$  
such that $\rho\!\!\pmod{\Fm_R}=\rhobar$, 
$\left.\rho\right|_{G_v}\in\CC_v$ for all $v\in S\cup T$ and~$\det\rho=\eta$.

\smallskip

To describe tangential conditions on the above lifts, we need 
to fix some more notation. For $v$ a place of $\widetilde X$ and $M$ a
$G_v$-module, the pairing $M\times M^*\to\BF$ defined by evaluation is
obviously perfect. Tate localy duality says that the induced pairing
\begin{eqnarray}\label{TracePairing0}
H^1(G_v,M)\times H^1(G_v,M^*(1))\longto H^2(G_v,\BF(1))\cong\BF
\end{eqnarray}
is perfect as well, and one denotes for any $\BF$-submodule $L\subset
H^1(G_v,M)$ its annihilator under this pairing by $L^\perp\subset
H^1(G_v,M^*(1))$. In the particular case of the subspace of unramified
cocycles 
$$H^1_\unr(G_v,M):=H^1(G_v/I_v,M^{I_v})\subset H^1(G_v,M)$$ 
one finds $H^1_\unr(G_v,M)^\perp = H^1_\unr(G_v,M^*(1)).$ 

\smallskip

We now specialize to the situations of interest to us, i.e., to
$M=\Ad$ and $M=\ad$. Note first that $\Ad$ is self-dual via the
perfect trace pairing $\Ad\times\Ad\to\BF:(A,B)\mapsto \Trace(AB)$. 
Because $\ell$ does not divide $n$, this pairing restricts to a
perfect pairing
\begin{eqnarray}\label{TracePairing}
\ad\times\ad\longto \BF,
\end{eqnarray}
and local Tate duality induces the perfect pairing
\begin{eqnarray}\label{TracePairing2}
H^1(G_v,\ad)\times H^1(G_v,\ad(1))\longto H^2(G_v,\BF(1))\cong\BF.
\end{eqnarray}

For a finite subset $T$ of $X$ and a collection 
$(L_v)_{v\in S\cup T}$ of subspaces of
$H^1(G_v,\ad)$ one defines $H^1_{\{L_v\}}(T,\ad)$ as the kernel of
$$H^1(\pi_1(X\setminus T),\ad)\longto \oplus_{v\in S\cup T}H^1(G_v,\ad)/L_v.$$
Note that $H^1_{\{L_v\}}(T,\ad)$ does depend on $S$ even though $S$
does not explicitly appear in the notation.

Ramakrishna's first observation is the following:
\begin{lemma}\label{FirstObs}
Suppose one is given locally admissible pairs 
$(\CC_v,L_v)_{v\in S\cup T}$ compatible with $\eta$ such that
$$H^1_{\{L_v^\perp\}}(T,\ad(1))=0.$$ 
Then there exists a lift of $\rhobar$ to $W(\BF)$ of type
$(\CC_v)_{v\in S\cup T}$. 
\end{lemma}

The proof is essentially that of \cite{taylor}, Lemma~1.2, and so we omit the
details.

\begin{remark}\label{WilesFormula}
{\em Mimicking the proofs of \cite{ddt}, Thms.~2.13, 2.14, one obtains
  for a $\pi_1(X\setminus T)$-module $M$ and subspaces $L_v\subset
  H^1(G_v,M)$ for $v\in S\cup T$ the formula} 
$$\frac{|H^1_{\{L_v\}}(T,M)|}{|H^1_{\{L_v^\perp\}}
(T,M^*(1))|}=\frac{|H^0(\pi_1(X),M)|}{|H^0(\pi_1(X),M^*(1))|}
\prod_{v\in S\cup T}
\frac{|L_v|}{|H^0(G_v,M)|}.
$$

{\em In our situation $M\cong  M^*\cong\ad$, the first quotient on the
  right is clearly $1$. Thus by Remark~\ref{Heuristically}, one
expects the product on the right to have the value at most
$1$. Furthermore, this should happen precisely when $\dim L_v=
h^0(G_v,\ad)$ for all $v\in S\cup T$. Therefore if the hypothesis of the
above lemma are satisfied, then one expects 
$$\dim H^1_{\{L_v\}}(T,\ad)=0.$$  
In terms of deformation theory, this can be 
interpreted by saying that the universal deformation ring of type
$(\CC_v)_{v\in S\cup T}$ is smooth over $W(\BF)$ of relative
dimension zero, i.e., isomorphic to $W(\BF)$. 

Note that the above formula also holds for $S\cup T=\emptyset$, even though
the duality results in \cite{nsw} are not proved in this case. 
The reason is that in this case the right hand side is $1$
and because of $H^0(\pi_1(X),\ad)=0$, the left hand expresses that 
fact that the Euler-Poincar\'e characteristic of the unramified 
$\BF[\pi_1(\widetilde X)]$-module $\ad$ is zero. }
\end{remark}

We need to slightly generalize the concept of sufficiently many
$R$-classes for the following result: Suppose we are given locally
admissible  $(\CC_v,L_v)_{v\in S\cup T}$ which are compatible with $\eta$.
We say that $\rhobar$ {\em admits sufficiently many $R$-classes for
  $(\CC_v,L_v)_{v\in S\cup T}$} if there exists at least one $R$-class and
if the kernels of the restriction morphisms (composed with $\ad\to
\ad_\sigma$ at $R$-places) 
$$H^1(\Gal(E(\zeta_\ell)/K),\ad)\cap H^1_{\{L_v\}}(T,\ad)\to\!\!\!\!\!\!
\prod_{\genfrac{}{}{0pt}{}{\sigma\,\mathrm{an}\,R\mathrm{-class}}{\mathrm{of}\,\mathrm{type(II)}}}
\!\!\!\!\! H^1(\langle\sigma\rangle,\ad_\sigma),$$
$$H^1(\Gal(E(\zeta_\ell)/K),\ad(1))\cap
H^1_{\{L^\perp_v\}}(T,\ad(1))\to\!\!\!\!\!\!
\prod_{\genfrac{}{}{0pt}{}{\sigma\,\mathrm{an}\,R\mathrm{-class}}{\mathrm{of}\,\mathrm{type(II)}}}
\!\!\!\!\! H^1(\langle\sigma\rangle,\ad_\sigma)$$
are zero.

The main observation of Ramakrishna, if adapted to our situation,
is the following key lemma:
\begin{lemma}\label{keylemma}
Suppose one is given a finite set of
places $T'\subset X$ and locally admissible  $(\CC_v,L_v)_{v\in S\cup T'}$
which are compatible with $\eta$ and such that  
$$\sum_{v\in S}\dim L_v\geq \sum_{v\in S}h^0(G_v,\ad).$$
If $\rhobar$ admits sufficiently many $R$-classes for
$(\CC_v,L_v)_{v\in S\cup T'}$, 
then one can find a finite set of $R$-places $T\subset X$ and
locally admissible $(\CC_v,L_v)_{v\in T}$ compatible with $\eta$, such
that $$H^1_{\{L_v^\perp\}}(T\cup T',\ad)=0.$$
\end{lemma}
The proof of this lemma will be given in Section~\ref{KeyLemSec}. Let
us now explain how this will give a proof of Theorem~\ref{OnSerreConj}:

In the following two sections, we will define good local
deformation problems at certain unramified primes and at ramified
primes where the ramification is either of order prime to $\ell$ or
prime to $p$. We then apply Lemma~\ref{keylemma} with $T'=\emptyset$
and assume that $\rhobar$ ramifies at all places of $S$. In order to
do that we also 
have to check that if $\rhobar$ has sufficiently many $R$-classes,
then this implies that $\rhobar$ has sufficiently many $R$-classes for
$(\CC_v,L_v)_{v\in S}$ where the $(\CC_v,L_v)$ will be defined
below. Once this is shown the theorem follows easily from the above
two lemmas. The full proof is given at the end of
Section~\ref{KeyLemSec}. 

\subsection{Local deformations at $R$-places}\label{AuxPrimes}
In this section, we will define locally admissible
deformation problems $\CC_v$ at $R$-places $v$ compatible with 
the Teichm\"uller lift $\eta_v\!:G_v\to W(\BF) $ of
$\det\rhobar_{v}$, cf.\ Proposition~\ref{MainPropOnAuxPrimes}. 
So for the remainder of this section, we fix an $R$-place $v$, 
and denote by $\sigma$ the image of
$\Frob_v$ in $\Gal(E(\zeta_\ell)/K)$, so that $[\sigma]$ is an $R$-class.
We also fix an eigenvalue $\lambda\in\BF$ of $A:=\rhobar(\sigma)$
as required in the definition of an $R$-place.

\subsubsection*{The definition of $\CC_v$ at an $R$-place}
Using the rational canonical form, we may assum that $A$ is given in
the form
$$ A=\left(\begin{array}{ccc}A_1&&0\\
&\ddots&\\ 0&&A_r \end{array}\right),$$
where each $A_i$ is a square matrix of size $n_i$, the matrices $A_i$
for $i>1$ are in rational canonical form and act indecomposably, and
the matrix $A_1$ has the following form depending on our our two
cases:
$$ A_1=\left(\begin{array}{cc}\lambda\chi(\sigma)&0\\
0&\lambda \end{array}\right)\ \ \hbox{ in case (I) }\quad\hbox{and}\quad 
A_1=\left(\begin{array}{cc}\lambda&1\\
0&\lambda
\end{array}\right)\ \ \hbox{ in case (II) }.$$
Note that in case (II) the $A_i$ are in bijection with the irreducible
representations $V_i$ used in the defining formula (\ref{DefOfAdSig}).
Because the $A_i$ act indecomposably, the eigenvalues form a single
Galois orbit and the Jordan canoncal form of an $A_i$ consists of
identical blocks for each of the eigenvalues. Because $A$ is regular,
different $A_i$ have distinct orbits of eigenvalues. Also clearly each
$A_i$ is again regular. 

For $i=2,\ldots,r$ we define $$\rho_{v,i}\!:G_v\longto \GL_n(R_{v,i})$$
as the versal unramified deformation of $\rhobar_{v,i}\!:G_v\to
\GL_n(\BF)$ defined as the restriction of $\rhobar$ to the $i$-th
block. 

For the definition in case $i=1$,
let $\hat\BZ$ be the profinite completion of $\BZ$ and
$\hat\BZ'$ the prime-to-$p$ completion of $\hat\BZ$. Let $s,t$ be
topological generators of $\hat\BZ$ and $\hat\BZ'$, respectively.
For $q$ a power of $p$, and thus prime to
$\ell$, define $\overline G_q:=\hat\BZ'\rtimes\hat\BZ$, where the
semidirect product is given (in multiplicative notation) by the
condition $sts^{-1}=t^q$. Then $\overline G_{q_v}$ can be identified
with the tame quotient of $G_v$ in such a way that $t$ is a generator
of $I_v/P_v$ and $s$ is a lift of $\Frob_v\in G_v/I_v$.

By $\hat\mu\in W(\BF)$ we denote the Teichm\"uller lift of any element 
$\mu$ of $\BF$, by $\tilde q_v:=q_v/\widehat{\chi(\sigma)}$, and we
set $\delta$ to be $0$ in case (I) and $1$ in case (II). We now define
$R_{v,1}:=W(\BF)[x_{1,0},x_{1,1}]$ and  
$\rho_{v,1}\!:G_v\onto \overline G_{q_v} \longto \GL_n(R_{v,i})$ 
by 
$$ s\mapsto \left(\begin{array}{cc}\widehat{(\lambda\chi(\sigma))}
  \tilde q_v^{1/2}(1-x_1)&\delta\\
0&\hat\lambda \tilde q_v^{-1/2}(1-x_1)\end{array}\right)\quad\hbox{and}\quad 
t\mapsto \left(\begin{array}{cc}1&x_0\\
0&1\end{array}\right).$$
The necessary condition $\rho_{v,1}(s)\rho_{v,1}(t)=
\rho_{v,1}(t)^{q_v}\rho_{v,1}(s)$ can easily be verified.

Combining the above, we define 
$$R_v:=\mathop{\widehat\otimes}\limits_{i=1}^r
R_{v,i}\Big/\Big( \prod_{i=1}^r \det\rho_{v,i}(s)-\eta_v(s)\Big),$$
with $\hat\otimes$ formed over $W(\BF)$, and the corresponding
representation $\rho_v\!:G_v\longto \GL_n(R_v)$ as $\oplus \rho_{v,i}$
(where the entries are taken modulo the ideal generated by
$\prod_{i=1}^r \det\rho_{v,i}(s)-\eta_v(s)$).

\medskip

To investigate the resulting representations, we first need a simple
result on the individual $\rho_{v,i}$. For this we denote by $\Ad_i$
the adjoint representations of the $A_i$ and by the $\ad_i$ its
subrepresentation on trace zero matrices; i.e., in
case (II) they agree with those defined in (\ref{DefOfAdSig}). Then
obstruction theory easily shows the following:
\begin{lemma}\label{LocDefsLem}
Let $i$ be in $2,\ldots,r$. Then the versal deformation $\rho_{v,i}$ is smooth
over $W(\BF)$ of dimension $h^1(G_v,\Ad_i)=n_i$. If $\ell\notdiv n_i$
and if $\eta_i$ is any lift of $\det\rhobar_i$ to $W(\BF)$, then the
versal deformation of determinant equal to $\eta_i$ is smooth of
dimension $h^1(G_v,\ad_i)= n_i-1$.
\end{lemma}

\begin{cor}\label{RvIsSmooth} Assume that there exists an $i$ such
  that $\ell$ does not divide $n_i$. Then $R_v$ is
  smooth over $W(\BF)$ of relative dimension $n-1=h^0(G_v,\ad)$.
\end{cor}
Note that for $\ell\notdiv n$ such an $n_i$ always exists.
\begin{Proof}
By the preceeding lemma, the ring $\widehat\otimes_{i=1}^r R_{v,i}$ is
smooth of dimension $n$. Because the $A_i$ have distinct sets of
eigenvalues for different $i$ we have $h^1(G_v,\Ad)=\sum_i
h^1(G_v,\Ad_i)=n$. Moreover if one of the $n_i$ is not divisible by
$\ell$, then it is easy to see that $h^1(G_v,\ad)=h^1(G_v,\Ad)-1$. Let
$i_0$ be the corresponding index. 

We now prove the smoothness of $R_v$. The previous lemma applied to
$i_0$ say that there is a system of local coordinates of $R_{v,i}$
such that $\det\rho_{v,i_0}(s)=\eta_{i_0}(s)(1+x)$ where $x$ is one of
these coordinates. If we regroup the defining relation of $R_v$, it
yields therefore the relation
$$\eta_v(s)\prod_{i\neq i_0} \det\rho_{v,i}^{-1}(s)=\eta_{i_0}(s)(1+x),$$
and the variable $x$ does not occur on the left hand side. Thus the
relation eliminates the variable $x$, which is 
one of the local coordinates in a suitable set of such for the ring
$\hat\otimes_iR_{v,i}$. Because $\hat\otimes_iR_{v,i}$ is smooth over
$W(\BF)$ of relative dimension $n$, so is $R_v$ of relative
dimension~$n-1$.
\end{Proof}

\medskip

The following defines a pair $(\CC_v,L_v)$ compatible with
$\eta_v$: The functor $\CC_v\!:\CA\longto\mathbf{Sets}$ is given by
\begin{eqnarray*} R&\mapsto& \CC_v(R):=\left\{\rho\!:G_v\to\GL_n(R)
\bigm| \exists\alpha\in\Hom_\CA(R_v,R),\right.\\
&&\left.\phantom{\CC_v(R):=\left\{\rho\!:\right.} \exists
  M\in1+M_n(\Fm_R):\rho=M(\alpha\circ\rho_v)M^{-1}\right\}. 
\end{eqnarray*}
Moreover, if $\rho_0\!:G_v\to\GL_n(\BF[\eps]/(\eps^2))$ denotes the
trivial lift of $\rhobar$, the subspace $L_v\subset H^1(G_v,\ad)$ is
the set of $1$-cocycles  
$$\Big\{c:G_v\to\ad: g\mapsto \frac1\eps(\rho(g)\rho_0^{-1}(g)-I)
\,\Big| \; \rho\in\CC_v(\BF[\eps]/(\eps^2)) \Big\}$$
and $L_{v,\unr}\subset H^1_\unr(G_v,\ad)$ is the intersection $L_v\cap
H^1_\unr(G_v,\ad)$. 

\begin{prop}\label{LemOnLv}
\label{MainPropOnAuxPrimes}
\begin{enumerate}
\item \label{MPOAPi} $\dim_\BF L_v=1+\dim_\BF L_{v,\unr} =n-1$.
\item \label{MPOAPii} The pair $(\CC_v,L_v)$ satisfies conditions P1--P7
  of~\cite{taylor}
\end{enumerate}
\end{prop}

\subsubsection*{Verifying the axioms for $\CC_v$}
Our first aim is to prove the assertions on the dimensions in the
above proposition. For this we need a closer analysis of
$H^1(G_v,\ad)$. Observe first that by repeatedly applying the
Leray-Serre spectral sequence to $G_v\supset I_v\supset P_v$ and
$\ad$, one obtains the short exact sequence
\begin{equation}\label{SESforH1}
0\to \ad/(\sigma-1)\ad \to H^1(G_v,\ad)\to (\ad(-1))^\sigma\to 0
\end{equation}
and an isomorphism $H^i(G_v,\ad)\cong H^i(\overline G_{q_v},\ad)$.

The above short exact sequence can be interpreted
in terms of $1$-cocycles representing cohomology classes. Namely any
$1$ cocycle $c$ of $\overline G_{q_v}$ with values in $\ad$ is
uniquely determined by its values $c(s)$, $c(t)$. These are subject to
the conditions $c(s)\in\ad$ and $c(t)\in(\ad(-1))^s$, i.e.,
$c(t)\in\ad$ satisfies $s c(t)=\frac1q c(t)$. Furthermore the
$1$-coboundaries are precisely the $1$-cocycles with $c(s)\in
(s-1)\ad$ and $c(t)=0$.

To determine $L_v$, we state the following simple lemma without proof:
\begin{lemma}\label{LemOnH1Unr} 
Suppose $A\in\GL_n(\BF)$ is regular. Let $V$ be the $\BF[A]$-module on
$\BF^n$ defined by $A$ and suppose $V=\oplus V_i$ is a direct sum
decomposition as a $\BF[A]$-module. Let the corresponding $\Ad_i$ be
defined as above. Then $\Ad/[\Ad,A]=\oplus_i \Ad_i/[\Ad_i,A]$. 
\end{lemma}

\bigskip

\begin{Proofof}{Proposition~\ref{MainPropOnAuxPrimes}\ref{MPOAPi}}
Let us fix local coordinates $x_{i,j}$, $j=1,\ldots,n_i$ of the rings
$R_{v,i}$, $i=2,\ldots,r$. We also enumerate them, so that the
variable $t$ in the proof of Corollary~\ref{RvIsSmooth} is given by
$x_{i_0,1}$.

Let first $c_0$ be the $1$-cocycle that arises from $R_v\to
\BF[\eps]/(\eps^2)$ by mapping the $x_{i,j}$, $j\ge1$ and $(i,j)\neq
(i_0,1)$, to zero and $x_{1,0}$ to $\eps$. (The image of $x_{i_0,1}$
is determined by the $x_{i,j}$ with $j\ge1$.) The correspondonding
element in $L_v$ is easily seen to non-zero and ramified. Moreover all
cocycles obtained from an assignment where $x_{1,0}$ maps to zero are
unramified. This shows $L_v=L_{v,\unr}+1$.

Let $L_v'$ be the set of cocycles corresponding to the ring
$R':=\hat\otimes_iR_{v,i}/(x_{1,0})$. The subspace defined from 
$R_{v,1}/(x_{1,0})$ is $1$-dimensional, that from $R_{v,i}$, $i>1$, 
has dimension $n_i$ by Lemma~\ref{LocDefsLem}. 
The previous lemma implies $\dim_\BF L'=1+\sum_{i>1}
n_i=n-1$ which in turn is the relative dimension of $R'$ over
$W(\BF)$. Because $R_v/(x_{1,0})$ is a smooth quotient of $R'$ of relative
dimension $n-2$ this yields $\dim_\BF L_{v,\unr}=n-2$, as asserted.
\end{Proofof}

\begin{Proofof}{Proposition~\ref{MainPropOnAuxPrimes}\ref{MPOAPii}} 
The only non-trivial condition to verify is P4. For this let us first
prove the following:  
\begin{lemma}
Let $\widetilde R$ be in $\CA$ and $\alpha,\alpha'\in\Hom_\CA(R_v,\widetilde
R)$ such that there exists $M\in\GL_n(\widetilde R)$ congruent to the
identity modulo $\Fm_{\widetilde R}$ with
$M(\alpha\circ\rho_v)M^{-1}=\alpha'\circ\rho_v.$
Then $\alpha\circ\rho_v(s)=\alpha'\circ\rho_v(s)$, so that in
particular $M$ commutes with $\alpha\circ\rho_v(s)$.
\end{lemma}
\begin{Proof}
We use the same local parameters for $R_v$ as in the proof of
Proposition~\ref{MainPropOnAuxPrimes}~(i). The matrix $\rho(s)$ has
entries in the power series ring over $W(\BF)$ in the variables
$x_{i,j}$. By $\rho(s)_r$ we denote the part of $\rho(s)$ that is
homogeneous of degree $r$, so that $\rho(s)=\sum_{r=0}^\infty
\rho(s)_r$. The assertion $\dim_\BF L_{v,\unr}=n-2$ of 
Proposition~\ref{MainPropOnAuxPrimes}~(i) means precisely
that the $n-2$ matrices $\frac{\partial}{\partial x_{i,j}}\rho(s)_1$
over all $i,j$ with $j\ge1$ and $(i,j)\neq(i_0,1)$ form a basis
of the vector space~$L_{v,\unr}\subset \ad/(s-1)\ad$.

Define $\rho:=\alpha\circ\rho_v$, $\rho':=\alpha'\circ\rho_v$. 
Let $\rho_{(m)}:=\rho\!\!\pmod{\Fm_{\widetilde R}^m}$ and introduce
analogous abbreviations $\rho'_{(m)}$, $\alpha_{(m)}$,
$\alpha'_{(m)}$ and $M_{(m)}$. By induction on $m$ we will show that
$\alpha_{(m)}(x_{i,j})=\alpha'_{(m)}(x_{i,j})$ for all $i,j$ with
$j\ge1$.  This clearly implies $\rho_{(m)}(s)=\rho'_{(m)}(s)$
for all~$i$, and thus the lemma. The case $m=1$ is clear and so we
now carry out the induction step~$m\mapsto m+1$.

By the
induction hypothesis $M_{(m)}$ commutes with $\rho_{(m)}(s)$. Because
$\rhobar(s)$ is regular, \cite{boeckle1}, Lemma 5.6, implies that
there exists a lift $M'$ of $M_{(m)}$ to $\GL_n(\widetilde
R/{\Fm_{\widetilde R}^{m+1}})$ which commutes with $\rho_{(m+1)}(s)$. 
By considering ${M'}^{-1}M_{(m+1)}$, we may thus assume $M_{(m+1)}=
I+\Delta$ for some $\Delta\in M_n(\Fm_{\widetilde R}^{m}
/\Fm_{\widetilde R}^{(m+1)})$. We also define elements 
$\delta_{i,j}:=\alpha_{(m+1)}(x_{i,j})-\alpha'_{(m+1)}(x_{i,j})$ which
by the induction hypothesis lie in $\Fm_{\widetilde R}^{m+1}$. The
expansion of $\rho(s)$ in homogeneous parts shows 
$$\rho_{m+1}'(s)=\rho_{m+1}(s)+\rho(s)_{1}|_{x_{i,j}=\delta_{i,j}}
\hbox{ in }M_n({\widetilde R}
/\Fm_{\widetilde R}^{(m+1)}),$$
and so the condition 
$M_{(m+1)}\rho_{(m+1)}(s)M^{-1}_{(i+1)}=\rho'_{(m+1)}(s)$ yields 
\begin{equation}\label{RecEqn}
\sum_{(i,j)} \delta_{i,j}\frac{\partial}{\partial x_{i,j}}\rho(s)_1=
\rho(s)_{1}|_{x_{i,j}=\delta_{i,j}}=
\Delta \rho_{(m+1)}(s)-\rho_{(m+1)}(s)\Delta 
\end{equation}
in $M_n(\Fm_{\widetilde R}^{m}/\Fm_{\widetilde R}^{(m+1)})$. The
right hand side is clearly a linear combination of coboundaries, of
$H^1(G_v,\ad)\otimes \Fm_{\widetilde R}^{m}/\Fm_{\widetilde
  R}^{(m+1)}$, while the left hand side is a linear combination of a
basis of the cohomology. Therefore both sides must vanish. This
concludes the induction step.
\end{Proof}
To verify P4, suppose that 
we are given rings $R_1,R_2\in\CA$, lifts $\rho_i\in\CC_v(R_i)$,
ideals $I_j\in R_j$, and an identification $R_1/I_1\cong R_2/I_2$
under which $\rho_1\!\!\!\!\pmod{I_1}\equiv \rho_2\pmod{I_2}$. We want to glue
the $\rho_i$ to an element $\rho$ of $\CC_v(R)$ for
$$R:=\{(r_1,r_2)\in R_1\oplus R_2:r_1\!\!\pmod{I_1}=r_2\!\!\pmod{I_2}\}.$$ 
So let $\alpha_i\in\Hom_\CA(R_v,R_i)$ and $M_i\in\GL_n(R_i)$ such that
$\rho_i=M_i(\alpha_i\circ\rho_v)M_i^{-1}$, $i=1,2$. We claim that
there exists $\alpha\in\Hom_\CA(R_v,R)$ and $M\in\GL_n(R)$ with
$M\equiv I\!\!\pmod {\Fm_R}$ such that
$\rho:=M(\alpha\circ\rho_v)M^{-1}=\rho_1\oplus\rho_2$. 

By conjugating $\rho_1$ by some lift of $M_2\!\!\pmod {I_1} $ to
$R_1$, we may assume that $M_2=I$. 
By the lemma, the matrix $M_1\!\!\pmod {I_1}$ 
commutes with $(\alpha_1\!\!\pmod {I_1})\circ\rho_v(s)=(\alpha_2\!\!\pmod
{I_2})\circ\rho_v(s)$. Using \cite{boeckle1}, Lemma 5.6, and the regularity of
$\rhobar$, we may choose a lift $M'_1\in\GL_n(R_1)$ of 
$M_1\!\!\pmod {I_1}$ which commutes with $\alpha_1\circ\rho_v(s)$. We
now replace $M_1$ by $\widetilde M_1:=M_1{M_1'}^{-1}$ and $\alpha_1$ by
some $\tilde\alpha_1\!:R_v\to R_1$, which differs from $\alpha_1$ at
most on the variable $x_0$, and such that 
$$ \widetilde M_1(\tilde\alpha_1\circ\rho_v)\widetilde M_1^{-1}=
M_1(\alpha_1\circ\rho_v)M_1^{-1}.$$
Defining $M:=(\widetilde M_1,I)\in\GL_n(R)$ and
$\alpha:=(\alpha_1',\alpha_2)\!:R_v\to R$, the above claim is
satisfied, and the proof of \cite{taylor}, P4, completed. 
\end{Proofof}

\subsubsection*{On the local duality pairing}

In analogy to the short exact sequence (\ref{SESforH1}) one also has the
short exact sequence
$$
0\to \ad(1)/(\sigma-1)\ad(1) \to H^1(G_v,\ad(1))\to \ad^\sigma\to 0.
$$
Thus one may view $L_{v,\unr}$ as a subspace of $\ad/(\sigma-1)\ad $
where $\sigma$ is the image of $\Frob_v$ in $\Gal(E(\zeta_\ell)/K)$
and $L_{v,\unr}^\perp:=L_v^\perp\cap H^1_\unr(G_v,\ad(1))$ as a
subspace of $\ad(1)/(\sigma-1)\ad(1)$. While it is clear that the
former only depends on $\sigma$ and the choice of $\lambda$, 
for the latter this is not
immediate, since it was defined using the pairing~(\ref{TracePairing2}) 
which in turn was defined using Tate local duality. For later
use, we now show that in fact also $L_{v,\unr}^\perp\subset
\ad(1)/(\sigma-1)\ad(1)$ only depends on $\sigma$. For this, let us
fix $v$ and $\sigma$ as above and also $\lambda$ and, in case
(I),~$\lambda'$.  

Above we observed that cocycles in $H^i(G_v,\ad)\cong H^i(\overline
G_{q_v},\ad)$ are determined by their values $c(s)\in \ad/((s-1)\ad)$ and
$c(t)\in(\ad(-1))^s$. Similar interpretations hold for $1$-cocycles in
the cohomology group $H^i(G_v,\ad(1)) \cong H^i(\overline G_{q_v},\ad(1))$. We will now
make explicit the Tate local duality in terms of cocycle
representatives. This requires that we make explicit the
pairing~(\ref{TracePairing2}), as well as the isomorphism
$H^2(G_v,\BF(1))\cong\BF$.

\smallskip

It is well-known and, using the Leray-Serre spectral sequence and the
fact that $\hat\BZ$ and $\hat\BZ'$ are of cohomological dimension one, easy to
see that $H^2(G_v,\BF(1))\cong
H^2(\overline G_q,\BF(1))\cong\BF$. Recall a) that elements of
$H^2(\overline G_q,\BF)$ classify extensions of
$\overline G_q$ by $\BF$ and b) that elements of this
cohomology group may be represented by normalized $2$ cocycles, i.e., 
maps $$[.,.]\!:\overline G_q\times \overline G_q\to\BF$$ which satisfy
$[1,g]=[g,1]=0$ for all $g\in\overline G_q$ and 
$$f[g,h]-[fg,h]+[f,gh]-[f,g]\quad\forall f,g,h\in \overline G_q.$$
Note that in our situation $f[g,h]=[g,h]$; but we decided to leave $f$
in the notation, to remind the reader of the condition of a normalized
$2$-cocycle also in the case of non-trivial coefficients. 

Regarding the duality pairing~(\ref{TracePairing2}), we have the
following results: 
\begin{lemma}\label{Lem1}
In terms of normalized $2$-cocycles, the isomorphism 
$$H^2(\overline G_q,\BF(1))\stackrel\cong\longto \BF$$
is given by 
$$[.,.]\mapsto
\sum_{i=1}^{q^{\ell-1}-1}[t^i,t]+[t^{q^{\ell-1}},s^{\ell-1}]-[s^{\ell-1},t]\in\BF.$$ 
\end{lemma}
\begin{lemma}\label{Lem2}
With respect to the isomorphism of the previous lemma, the trace
pairing 
$$H^1(G_v,\ad)\times H^1(G_v,\ad(1))\longto \BF.$$
is given explicitly in terms of $1$-cocycles as follows: Let $c_1$ and
$c_2$ be $1$-cocycles of $ H^1(\overline G_{q_v},\ad)$ and $ H^1(\overline
G_{q_v},\ad(1))$, 
respectively. Then the image of $(c_1,c_2)$ under the pairing is given
by $$\mathrm{Trace}(c_1(s)c_2(t)- c_2(s)c_1(t))\in\BF,$$
unless $\ell=2$ and $q\equiv3\pmod4$. In the latter case it is
$\mathrm{Trace}(c_1(s)c_2(t)+ c_2(s)c_1(t)+c_1(t)c_2(t))$.
\end{lemma}
\begin{cor}\label{CorOnSigmaAndV}
$L_{v,\unr}$ and $L_{v,\unr}^\perp$ only depend on~$\sigma$ and the
choice of eigenvalue $\lambda$.
\end{cor}
Therefore, in later sections we freely write
$L_{\sigma,\unr}$ and $L_{\sigma,\unr}^\perp$ for the respective
spaces, where $\sigma$ is the image of $\Frob_v$ in $\Gal(E(\zeta_\ell)/K)$.
If we also want to include the choice of $\lambda$ in the notation, we
write $L_{\sigma,\lambda,\unr}$ and $L_{\sigma,\lambda,\unr}^\perp$.

Let us first prove the corollary:
\begin{Proof}
On rings of characteristic $\ell$, the definition of $\CC_v$ depends only
on $\sigma$. This easily implies that $L_{v,\unr}\subset
\ad/(\sigma-1)\ad$, as well as
$M:=L_v+H^1_\unr(G_v,\ad)/H^1_\unr(G_v,\ad)\subset(\ad(-1))^\sigma$
only depend on $\sigma$. Let $c_1$ be any cocycle in $L_v$
whose image $\tilde c_1$ in $M\subset (\ad(-1))^\sigma$ is a
generator. \ Since $H^1_\unr(G_v,\ad)$ and
$H^1_\unr(G_v,\ad(1))$ are orthogonal under the trace pairing,
Lemma~\ref{Lem2} shows that $L_{v,\unr}^\perp$ is spanned by the set
of all unramified cocycles $c_2\!:\overline G_{q_v}\to \ad(1)$ which satisfy 
$\mathrm{Trace}(\tilde c_1 \cdot c_2(s))=0$. This gives a condition on
$c_2(s)\in\ad(1)/(\sigma-1)\ad(1)$ which only depends on $M$, i.e.\ on
$\sigma$. The corollary is thus proved. 
\end{Proof}

\begin{Proofof}{Lemma~\ref{Lem2}}
We assume that we have proved Lemma~\ref{Lem1}. In terms of
$1$-cocycles, the map $$H^1(\overline G_q,\ad)\times H^1(\overline G_q,\ad))\to
H^2(\overline G_q,\ad\otimes\ad(1))$$ is given by mapping a pair $(c_1,c_2)$
to the (normalized) $2$-cocycle defined by $[f,g]:=c_1(f)\otimes
c_2(g)$. If we compose this with the map on cohomology  induced from the
trace map 
$$\ad\otimes\ad(1)\to \BF(1):A\otimes B\mapsto
\mathrm{Trace}(AB),$$
we obtain the (normalized) $2$-cocycle defined by
$[f,g]:=\mathrm{Trace}(c_1(f)c_2(g))\in\BF(1)$. By Lemma~\ref{Lem1} it
follows that the pair $(c_1,c_2)$ is mapped to 
$$ \mathrm{Trace}\Big(\sum_{i=1}^{q^{\ell-1}-1}c_1(t^i)c_2(t)+
c_1(t^{q^{l-1}})c_2(s^{\ell-1})-c_1(s^{\ell-1})c_2(t)\Big)\in\BF.$$
Because $c_1$ restricted to $\hat\BZ'$ is a homomorphism, we have
$c_1(t^i)=ic_1(t)$. So the sum simplifies to 
$$c_1(t)c_2(t)\sum_{i=1}^{q^{\ell-1}-1}i=c_1(t)c_2(t)q^{\ell-1}(q^{\ell-1}-1)/2.$$
As $q^{\ell-1}\equiv1\!\!\pmod \ell$, this sum is zero unless $\ell=2$
and $q\equiv3\pmod4$. In the latter case it is $c_1(t)c_2(t)$. For
the same reason, the term $c_1(t^{q^{\ell-1}})=q^{\ell-1}c_1(t)=c_1(t)$. 

To complete the proof of the lemma, it now suffices to show that we
may replace $c_2(s^{\ell-1})$ by $-c_2(s)$ (and similarly 
$c_1(s^{\ell-1})$ by $-c_1(s)$). An easy calculation shows that 
$\mathrm{Trace}(c_1(t)(\sigma-1)c_2(s))=0$. Also we have
$c_2(s^{\ell-1})=(1+\sigma+\ldots+\sigma^{\ell-2})c_2(s)$.
Combining the previous two observations, we find
$$\mathrm{Trace}(c_1(t)c_2(s^{\ell-1}))=\mathrm{Trace}(c_1(t)(\ell-1)c_2(s))
=-\mathrm{Trace}(c_1(t)c_2(s)),$$  
as asserted. The argument for $c_1(s^{\ell-1})$ is analogous.
\end{Proofof}

\begin{Proofof}{Lemma~\ref{Lem1}}
By the Leray-Serre spectral sequence applied to
$$\hat\BZ'\rtimes(\ell-1)\hat\BZ\subset \hat\BZ'\rtimes\hat\BZ$$
and the module $\BF(1)$, we obtain an isomorphism
$$H^2(\hat\BZ'\rtimes\hat\BZ,\BF(1))\cong
(H^2(\hat\BZ'\rtimes(\ell-1)\hat\BZ,\BF)(1))^{\BZ/(\ell-1)}$$
given by restriction. The point is that the action of
$\hat\BZ'\times(\ell-1)\hat\BZ$ on $\BF(1)$ is trivial (the residue field
of the corresponding local Galois extension has order $q^{\ell-1}$, and
hence contains a primitive $\ell-1$-th root of unity).
Since both $H^2(\ldots)$ terms are isomorphic to
$\BF$, it is not necessary to take invariants on the right for there
to be an isomorphism. So it suffices to show that the identification
asserted in the lemma is given by first restricting normalized
$2$-cocycles and then giving an isomorphism
$H^2(\hat\BZ'\rtimes(\ell-1)\hat\BZ,\BF)\cong\BF$. 

For the latter we use the interpretation in terms of extension
classes, cf.\ \cite{weibel}, \S 6.6. So let $[.,.]$ be a normalized $2$-cocycle
of the latter module. Then the corresponding extension $G$ can be
described as the group whose underlying elements are pairs $(a,x)$, $a\in\BF$,
$x\in\hat\BZ'\rtimes\hat\BZ$ and whose composition law is given by
$(a,x)(b,y)=(a+x\cdot b+[x,y],xy)$. One easily verifies that $G$ is
split if and only if there exist $a,b\in\BF$ such that $\tilde
s:=(a,s^{\ell-1})$ and $\tilde t:=(b,t)$ satisfy 
$$\tilde s\tilde t=\tilde t^{q^{\ell-1}}\tilde s.$$
In terms of $2$-cocycles, $G$ being split is equivalent to $[.,.]$
being a $2$-coboun\-da\-ry. 
Using the composition law, one can compute both sides. Let us
denote the difference of the $\BF$-component by $d([.,.])$, so that
$$d([.,.])=\sum_{i=1}^{q^{\ell-1}-1}[t^i,t]+[t^{q^{\ell-1}},s^{\ell-1}]-[s^{\ell-1},t].$$ 
We have $d([.,.])=0$ if and only if $[.,.]$ is a $1$-coboundary. 
Furthermore $d$ is $\BF$-linear and thus it induces an isomorphism 
$$H^2(\hat\BZ'\rtimes(\ell-1)\hat\BZ,\BF)\cong \BF.$$
Given a $2$-cocycle for $\hat\BZ'\rtimes\hat\BZ$, restricting it
to $\hat\BZ'\rtimes(\ell-1)\hat\BZ$ and applying $d$ yields precisely the
formula in the lemma, and so its proof is completed. 
\end{Proofof}

\subsection{Local deformations at $r$-places}\label{RamPrimes}

Regarding places at which $\rhobar$ is
ramified, one has the following results: 
\begin{prop}\label{FirstPropOnRamV}
Suppose that $\rhobar(I_v)$ is of order prime to $\ell$. Define 
the functor $\CC_v\!:\CA\to\mathbf{Sets}$ by 
$$R\mapsto \left\{\rho\!:G_v\to\GL_n(R)\;|\;\rho\!\!\!
\pmod{\Fm_R}=\rhobar_{v},\, 
\rho(I_v)\cong\rhobar(I_v),\,\det\rho=\eta_v \right\}$$
and $L_v$ as the corresponding subspace in $H^1(G_v,\ad)$. 
Then $(\CC_v,L_v)$ satisfies the conditions P1--P7 of \cite{taylor}, 
the conductors of $\rhobar_{v}$ and of any lift $\rho\in\CC_v(R)$,
$R\in\CA$ agree, $L_v=H^1_\unr(G_v,\ad)$ and $\dim L_v=h^0(G_v,\ad)$.
\end{prop}
\begin{Proof}
Except for the assertion on conductors, this is essentially
\cite{taylor}, Example E1, and so we only prove the latter part. For any ring
$R\in\CA$ and $\rho\in\CC_v(R)$, let $V_\rho(R)$ denote the $R[G_v]$
module defined by $\rho$. The kernel of $\GL_n(R)\to\GL_n(\BF)$ is a
pro-$\ell$ group and thus prime to $p$. Therefore the Swan conductors of
$\rho$ and $\rhobar_{v}$ are the same. The module $V_\rho(R)$ is free
over $R$ and thus the difference of the conductors of the two
representations is given by
$$\rank_R V_\rho(R)^{I_v}-\dim_{\BF}V_{\rhobar_{v}}(\BF)^{I_v}.$$
Since $G_v$ acts on both representations via the same quotient $\overline I_v$
which is prime to $\ell$, there is a natural equivalence between $R[\overline
I_v]$-representations which are free and finite over $R$ and $\BF[\overline
I_v]$-representations given by reduction modulo $\Fm_R$. In
particular both categories are semisimple and thus the above
expression is well-defined. Furthermore, this implies that the number
of trivial components contained in $V_\rho(R)$, as an $\overline
I_v$-module, is the same as that of $V_{\rhobar_{v}}(\BF)$, and
hence that the above difference is zero, as asserted.
\end{Proof}

\begin{prop}\label{SecondPropOnRamV}
Suppose that $\rhobar_v$ is at most tamely ramified and that
$h^0(G_v,\ad)<h^0(G_v,\Ad)$. Then there exists a pair
$(\CC_v,L_v)$ which satisfies conditions P1--P7 of \cite{taylor} with 
$L_v=H^1_\unr(G_v,\ad)$ and $\dim L_v =h^0(G_v,\ad)$, is compatible
with~$\eta_v$ and such that the conductors of $\rhobar_{v}$ and any
lift $\rho\in\CC_v(R)$, $R\in\CA$, agree.
\end{prop}

\begin{remark} One can construct examples which show that the
  condition $\dim_\BF\ad^{G_v}<\dim_\BF\Ad^{G_v}$ is necessary. The
  latter is automaticially satisfied if~$\ell\notdiv n$.
\end{remark}

\begin{remark}\label{RemOnTwoLvCv}{\em
If $\rhobar(I_v)$ is of order prime to $\ell p$, we may
apply either Proposition~\ref{FirstPropOnRamV} or
Proposition~\ref{SecondPropOnRamV} to obtain a pair $(\CC_v,L_v)$. The
pairs so obtained do have similar properties and in fact, in 
Remark~\ref{DefsAreEq} we will explain why the two are isomorphic.
}
\end{remark}

\medskip

In the remainder of this section, we shall give the proof of
Proposition~\ref{SecondPropOnRamV}. Since all representations that
occur in the proof will factor via the tame quotient $\overline
G_{q_v}$ of $G_v$, we fix the usual generaters $s$, $t$ satisfying the
relation $sts^{-1}=t^{q_v}$.

For $B\in \GL_n(W(\BF))$ we denote by $V$ the corresponding
$W(\BF)[X]$-module on $W(\BF)^n$ by having $X$ act via $B$.
Let $\BQ_\BF$ denote the fraction field of $W(\BF)$. 
We say that $B\in \GL_n(W(\BF))$ is a {\em  minimal lift} of its
reduction $\overline B\in  \GL_n(\BF)$, if $V=\oplus_i
V_{i,s}\otimes_{W(\BF)}V_{i,u}$ where the $V_{i,?}$ are
$W(\BF)[X]$-modules such that:
\begin{enumerate}
\item on $V_{i,u}$ the matrix representing $X$ is $W(\BF)$-conjugate
  to a regular unipotent matrix in Jordan form,
\item on $V_{i,s}$ the characteristic polynomial of $X$ is irreducible
  and its roots are Teichm\"uller lifts of elements in $\overline\BF$.
\end{enumerate}
\begin{lemma}
Any $\overline B\in\GL_n(\BF)$ has a minimal lift to $\GL_n(W(\BF))$.
\end{lemma}
\begin{Proof}
Let $\overline V:=\BF^n$ be the $\BF[X]$-module $\overline V:=\BF^n$
obtained by having $X$ acts as $\overline B$. We choose a
decomposition $\overline V\cong \oplus \overline V_i$ into
indecomposable $V_i$. On $V_i$ the action of $X$ decomposes into commuting semisimple and a
unipotent parts, defined over $\BF$. For instance by considering
Jordan normal forms over $\overline\BF$, one shows that
correspondingly one has $\overline V_i\cong \overline
V_{i,s}\otimes_\BF\overline V_{i,u}$ where $\overline V_{i,s}$ is a
semisimple representation of $X$ and $\overline V_{i,u}$ is a
unipotent representation of $X$. Because $\overline V_i$ is
indecomposable, the characteristic polynomial of $X$ on $\overline
V_{i,s}$ is irreducible over $\BF$. For the same reason, the action of
$X$ on $\overline V_{i,u}$ is by a regular unipotent matrix. So we may
assume that the operation of $X$ on $\overline V_{i,s}$ is given by a
companion matrix whose characteristic polynomial is irreducible over
$\BF$, and on $\overline V_{i,u}$ by a single Jordan block with
eigenvalue~$1$.

We now lift $X$ on $\overline V_{i,u}$ to a single Jordan block with
eigenvalue~$1$ over $W(\BF)$ and $X$ on $\overline V_{i,u}$ to a
companion matrix with eigenvalues the Teichm\"uller lifts of those of
$X$ on $\overline V_{i,u}$. The corresponding representations
$V_{i,u}$ and $V_{i,s}$ combine to give a representation of
$W(\BF)[X]$ on $V=\oplus_i V_{i,s}\otimes_{W(\BF)}V_{i,u}$ which has
all the required properties. Therefore the matrix representing this
$X$ is a minimal lift of~$\overline B$.
\end{Proof}

We first prove the following result, of which part (a) in the case
where $\rho(I_v)$ is an $\ell$-group is \cite{boeckle2}, Proposition~3.2.
\begin{prop}\label{MyOldProp}
Let $B$ be a minimal lift of $\overline B:=\rhobar(t)$.
\begin{enumerate}
\item $\CM:=M_n(W(\BF))/\{AB-B^{q_v}A\bigm|A\in
  W(\BF)\}$ is flat over~$W(\BF)$.
\item There exists a lift $\rho_0\!:G_v\onto\overline
G_v\longto\GL_n(W(\BF))$ with $\rho_0(t)=B$. 
\end{enumerate}
\end{prop}
\begin{Proof}
Let $X$ and the $V_{i,?}$ be as in the definition of minimal lift of
$B$. It is not  
difficult to see from condition~(ii) that we may assume that $X$ on
$V_{i,s}$ is given as a companion matrix $B_{i,s}$. Let now $\BF'$ be
a finite extension of $\BF$ which contains all eigenvalues of
$\overline B$. Then clearly over $W(\BF')$ the companion matrices
$B_{i,s}$ may be diagonalized. Moreover this diagonalization procedure
commutes with reduction modulo $\ell$. Since the base change
$\underline{\phantom{m}}\otimes_{W(\BF)}W(\BF')$
is faithfully flat, we will from now on for the proof of (a) assume
that $\BF$ contains all the eigenvalues of $\overline B$. 

Reordering
matrices in Jordan form and using the relation $\overline A\,\overline
B\,\overline A^{-1}=\overline B^{q_v}$ it is easy to prove the
following lemma. We leave details to the reader.
\begin{lemma}\label{NormalizationAtRamPrimes}
There exist $\mu_i\in\BF$ with Teichm\"uller lifts $\hat\mu_i$ and
$m_i\in\BN$, $i=1,\ldots,d$, such that \begin{enumerate}
\item \label{NormalizationAtRamPrimesPart1}$\mu_i^{q_v^{m_i}}=\mu_i$,
\item \label{NormalizationAtRamPrimesPart2}the elements $\mu_i^{q_v^j}$,
$i=1,\ldots,d$, $j=1,\ldots,m_i$ are pairwise disjoint and form a complete
list of the eigenvalues of $\overline B:=\rhobar(t)$, and 
\item with respect to a suitable basis one has
$$B=\left(\begin{array}{ccc}B_1&&0\\
&\ddots&\\ 0&&B_d
\end{array}\right), \hbox{ where each }B_i =\left(\begin{array}{ccc}B_{i,1}&&0\\
&\ddots&\\ 0&& B_{i,m_i}
\end{array}\right)
$$
is a square matrix, and for fixed $i$ the $B_{i,j}$ can be written as
$B_{i,j}=\hat\mu_i^{q_v^{j-1}}U_{i}$, for some unipotent $U_{i}$
in Jordan form, independent of $j$.
\end{enumerate}
Furthermore if $B$ is given as above, then $\overline A:=\rhobar(s)$
takes the form
$$\overline A=\left(\begin{array}{ccc}\overline A_1&&0\\
&\ddots&\\ 0&&\overline A_d
\end{array}\right)\hbox{ \ with }
\overline A_i =\left(\begin{array}{ccccc}0&\overline A_{i,1}&0&\cdots&0\\
&0&\overline A_{i,2}&\ddots&\\ \vdots&&\ddots&\ddots&0\\
0&&&0&\overline A_{i,m_i-1}\\
\overline A_{i,m_i}&0&\cdots&&0
\end{array}\right)$$
such that the $\overline A_{i,j}$ satisfy the relation $\overline
A_{i,j}U_i=U_i^{q_v}\overline A_{i,j}$.
\end{lemma}

We apply Proposition~\ref{MyOldProp}~(b) in the case
$\ell\notdiv\#\rhobar(I_v)$, covered by 
\cite{boeckle2}, to obtain matrices $A_{0,i,j}$ over $W(\BF)$ that satisfy
$A_{0,i,j}U_i=U_i^{q_v} A_{0,i,j}$ for all $i,j$ and whose reduction
modulo $\ell$ agree with $\overline A_{i,j}$. Let $A_0$ be composed from
the $A_{0,i,j}$ in the same way as $\overline A$ is from the
$\overline A_{i,j}$. Then $A_0$ is a lift to $W(\BF)$ of $\overline A$
such that $A_0BA_0^{-1}=B^{q_v}$. (This proves (b) only under the
further hypothesis that $\BF$ contains all eigenvalues of~$B$.)

We now consider the exact sequence 
\begin{equation}\label{SES-ForLifting}
0\to \CK\to M_n(W(\BF))\stackrel{A\mapsto AB-B^{q_v}A}\longto
M_n(W(\BF)) \to \CM\to 0
\end{equation}
where $\CK:=\{A\in W(\BF):AB=B^{q_v}A\}$. To complete (a), we need to show that the
generic rank of $\CM$ is the same as its special rank. To simplify the
problem we apply the isomorphism $M_n(W(\BF))\to M_n(W(\BF)):A\mapsto
A_0^{-1}A$ to the middle terms in~(\ref{SES-ForLifting}). This yields
the isomorphic exact sequence
\begin{equation*}
0\to \CK'\to M_n(W(\BF))\stackrel{A'\mapsto A'B-BA'}\longto
M_n(W(\BF)) \to \CM'\to 0,
\end{equation*}
with kernel $\CK'=\{A'\in M_n(W(\BF)):A'B=BA'\}\cong \CK'$,
and cokernel $\CM'=M_n(W(\BF))/\{A'B-BA'\bigm|A'\in W(\BF)\}\cong
\CM$.  We need to prove that the generic and special ranks of
$\CM'$ agree. 

Counting dimensions in the above short exact sequence and its
reduction modulo $\ell$, it suffices to show that the
dimension of $\CK'\otimes_{W(\BF)}\BQ_\BF$ and of $\overline\CK':=\{\overline
A'\in M_n(\BF)\bigm| \overline A'\,\overline B =\overline B\,\overline A'\}$ agree.
Because the $B_{i,j}$ have distinct eigenvalues modulo $\ell$, the
matrices $\overline A'\in\overline \CK'$ and $A'\in \CK'$, respectively,
will have the same block form as $B$. So we may consider blocks for
each pair $i,j$ separately. Therefore it suffices to prove the
assertion in the case where $B$ is a single Jordan block with
eigenvalue $1$. This case was treated explicitly in the proof of
\cite{boeckle2}, Proposition~3.2. The proof of (a) is now complete.

\smallskip

It remains to deduce (b) from (a). Because $\CM$ is flat, the
reduction mod $\ell$ of the exact sequence \ref{SES-ForLifting}
remains exact, and so the kernel of the reduction is $\CK/\ell\CK$.
The matrix $\overline A=\rhobar(s)\in M_n(\BF)$ lies in this kernel
and is therefore the reduction modulo $\ell$ of a matrix $A\in
\CK$. Because $A$ and $B$ satisfy the same relations as $s,t$ the
desired lift exists.
\end{Proof}

As a corollary to the above proof, we record the following technical
result, obtained by base change and using flatness.
\begin{cor}\label{CorOnSecondRamProp}
Suppose $B\in\GL_n(W(\BF))$ is a minimal lift of $\overline
B\in\GL_n(\BF)$. Then for any $R\in\CA$ the submodules
$\CK(R):=\{A\in M_n(R):AB=B^{q_v}A\}$ and $\CK'(R):=\{A\in
M_n(R):AB=BA\}$ of $M_n(R)$ are free and direct summands of $R$-rank
independent of $R$. Moreover $\CK'(R)=A_0^{-1}\CK(R)$ for $A_0$ as in
the above proof.
\end{cor}

\begin{Proofof}{Proposition~\ref{SecondPropOnRamV}}
Let $b_1,\ldots,b_m$ be parts of a basis of $\CK(W(\BF))$. We shall
specify more requirements on these elements below.
Let $x_1,\ldots,x_m$ be indeterminates and define
$$\CS_v:=A+\sum x_ib_i\in\GL_n(R), \CT_v:=B,$$
$$R_v:=W(\BF)[[x_1,\ldots,x_m]]/(\det\CS_v-\eta_v(s)),$$
$$\rho_v\!:G_v\onto\hat\BZ'\rtimes\hat\BZ\longto\GL_n(R_v):s\mapsto
\CS_v,t\mapsto \CT_v,$$
and the functor $\CC_v\!:\CA\to\mathbf{Sets}$ by 
\begin{eqnarray*} R&\mapsto& \CC_v(R):=\left\{\rho\!:G_v\to\GL_n(R)\,|\,
\exists\alpha\in\Hom_\CA(R_v,R),\right.\\
&&\left.\phantom{\CC_v(R):=\left\{\rho\!:\right.} 
\exists M\in1+M_n(\Fm_R):\rho=M(\alpha\circ\rho_v)M^{-1}\right\}.
\end{eqnarray*}

Let $L_v$ be the corresponding subspace $L_v\subset H^1(G_v,\ad)$. As
$\rho_v(t)$ does not deform, the subspace $L_v$ lies inside
$$H^1_\unr(G_v,\ad)\cong {\ad}^{t}/(s-1)\ad^{t}.$$ 
We denote by $\bar b_1,\ldots,\bar b_m\in\CK(\BF)$ the reductions of
the $b_i$ modulo $\ell$. So the elements $\overline A^{-1}\bar b_i$
lie in $\CK'(\BF)=\ad^t$, and an explicit calculation shows that 
$L_v$ is spanned by the images in ${\ad}^{t}/(s-1)\ad^{t}$ of those
linear combinations $\sum x_i\overline A^{-1}\bar b_i$, $x_i\in\BF$, which are
consistent with the determinant condition $\det\CS_v=\eta(s)$.

The latter condition modulo $(\ell,\Fm_R^2)$ means
\begin{eqnarray*}\det \overline A&\stackrel!=
&\det \overline A\cdot \det(I+\overline A^{-1}(\sum x_i\bar b_i)),\hbox{ i.e.,}\\
1&=&1+\sum x_i\mathrm{Trace}(\overline A^{-1}\bar b_i).
\end{eqnarray*}
Thus the above linear combinations satisfy $\sum
x_i\mathrm{Trace}(\overline A^{-1}\bar b_i)=0$. 

Now we fix the choice of the $b_i$. Namely, we take them as a subset
of $\CK(W(\BF))$ whose reductions modulo $\ell$ forms a basis of 
$\CK(\BF)/\{\overline X\,\overline A -\overline A\, \overline X
:\overline X\in\CK(\BF)\}.$ 
The elements $\overline A^{-1}\bar b_i$ then form a basis of
$$\Ad^t/(s-1)\Ad^{t}=\CK'(\BF)/\{\overline A \,\overline X\,\overline
A^{-1}-\overline X:\overline X\in\CK'(\BF)\}.$$
To pass from $\Ad$ to $\ad$ we use our assumption
$$h^0(G_v,\Ad)>h^0(G_v,\ad).$$ Because $(\Ad^{I_v})^{G_v/I_v}\cong
(\Ad^{I_v})_{G_v/I_v}$ as an $\BF$-vector space, and similarly for
$\ad$, we deduce that $(\Ad^t)/(s-1)\Ad^{t}$ properly contains
$(\ad^t)/(s-1)\ad^{t}$. This in turn shows that any element of 
the module $(\ad^t)/(s-1)\ad^{t}$ can be obtained as a linear combination $\sum
x_i\overline A^{-1}\bar b_i$ which satisfies $\sum_i
x_i\Trace(\overline A^{-1}b_i)=0$. This has two consequences:
Firstly we have $L_v=H^1_\unr(G_v,\ad)$; secondly the relation $\det
\CS_v=\eta(s)$ allows one to eliminate one of the variables $x_i$,
since this is possible tangentially, and so $R_v$ is smooth of
relative dimension $\dim L_v$ over~$W(\BF)$.

Note also that the determinant of $\rho_v$ is the Teichm\"uller lift
of that of $\rhobar_v$: For $\rho_v(t)$, this follows from the
construction of $B$, for $\rho_v(s)$ from the definition
of~$R_v$. This implies the result in general, since $\rho_v$ is only
tamely ramified.

\smallskip

Let us now verify properties P1--P7 of \cite{taylor}. As expected the
only property that is non-trivial is P4. To verify it, 
suppose we are given rings $R_1,R_2\in\CA$, lifts $\rho_i\in\CC_v(R_i)$,
ideals $I_j\in R_j$, and an identification $R_1/I_1\cong R_2/I_2$
under which $\rho_1\!\!\!\!\pmod{I_1}\equiv \rho_2\pmod{I_2}$. We need
to show that $(\rho_1,\rho_2)$ lies in $\CC_v(R)$ for
$$R:=\{(r_1,r_2)\in R_1\oplus
R_2:r_1\!\!\pmod{I_1}=r_2\!\!\pmod{I_2}\}.$$

So let $\alpha_i\in\Hom_\CA(R_v,R_i)$ and $M_i\in\GL_n(R_i)$ such that
$\rho_i=M_i(\alpha_i\circ\rho_v)M_i^{-1}$, $i=1,2$. We claim that
there exists $\alpha\in\Hom_\CA(R_v,R)$ and $M\in\GL_n(R)$ with
$M\equiv I\!\!\pmod {\Fm_R}$ such that
$(\rho_1,\rho_2)=M(\alpha\circ\rho_v)M^{-1}$. 
By conjugating $\rho_1$ by some lift of $M_2\!\!\pmod {I_1} $ to
$R_1$, we may assume~$M_2=I$. 

By an inductive argument, which is left to the reader, one can show
the following auxiliary result:
\begin{lemma}
Suppose $\widetilde R\in\CA$, $\widetilde J$ is a proper ideal of $\widetilde R$,
$A'\in A+\sum \beta_i b_i+M_n(\widetilde J)$ for some $\beta_i\in \widetilde
R$ and that $A'B=B^{q_v}A'$. Then there exists $\beta'_i\in \widetilde R$
with $\beta'_i-\beta_i\in \widetilde J$ and  $C\in I+M_n(\widetilde J)$ such
that $$A'=C\Big(A+\sum\beta_i' b_i\Big)C^{-1}.$$
\end{lemma}

Continuing with the proof of Proposition~\ref{SecondPropOnRamV},
observe that the condition 
\begin{equation}\label{EqnForP4}
M_1(\alpha_1\circ\rho_v)M_1^{-1}\!\!\!\pmod{I_1}=\alpha_2\circ\rho_v
\!\!\!\pmod{I_2}
\end{equation}
applied to $t$ implies that $M_1\!\!\pmod{I_1}$ commutes with
$B$. By Corollary~\ref{CorOnSecondRamProp}, we can find a lift
$\widetilde M_1\in M_n(R_1)$ of $M_1\!\!\pmod{I_1}$ which commutes
with~$B$. 

Because of (\ref{EqnForP4}) and the choice of $\widetilde M_1$, we can
apply the above lemma to $A':=\widetilde
M_1(\alpha_1\circ\rho_v(s))\widetilde M_1^{-1}$. It yields
$\alpha_1'\!:R_v\to R_1$ and $C\in I+M_n(I_1)$ such that 
$$C(\alpha_1'\circ\rho_v(s))C^{-1}=
\widetilde M_1(\alpha_1\circ\rho_v(s))\widetilde M_1^{-1}$$
and $\alpha_1'\!\!\pmod{I_1}=\alpha_2\!\!\pmod{I_2}$. Define
$\widetilde C:=M_1\widetilde M_1^{-1}C\in I+M_n(I_1)$. Then
$$\widetilde C(\alpha_1'\circ\rho_v)\widetilde C^{-1}=M_1(\alpha_1\circ\rho_v)M_1^{-1}=\rho_1$$
Therefore, if we set $M:=(\widetilde C,I)\in I+M_n(\Fm_R)$ and 
$\alpha:=(\alpha_1',\alpha_2)\!:R_v\to R$, we have
$$(\rho_1,\rho_2)=M(\alpha\circ\rho_v)M^{-1},$$
and so the proof of P4, and hence of all the axioms of Taylor, is completed.

\smallskip

It remains to prove the assertion on the conductors. As in the proof
of Proposition~\ref{FirstPropOnRamV}, the difference in conductors is
given by
$$\dim_R V_\rho(R)^{I_v}-\dim_{\BF}V_{\rhobar_{v}}(\BF)^{I_v},$$
where the notation is analogous to that in the quoted proof. Since
$I_v$ is topologically generated by the single element $t$, whose
image is the image of the matrix $B\in \GL_n(W(\BF))$, this difference
is given by 
$$ \dim_R \CK'(R)-\dim_{\BF}\CK'(\BF).$$
By Corollary~\ref{CorOnSecondRamProp} this difference is zero.
This shows that the conductors of $\rho$ and $\rhobar_v$ agree.
\end{Proofof}

\begin{remark}\label{DefsAreEq}
{\em Suppose now that the image of $I_v$ under $\rhobar$ is of order prime
to $\ell p$. Let $(\rho'_v,R'_v)$ be the versal deformation constructed
in Proposition~\ref{FirstPropOnRamV} and $(\rho_v,R_v)$ the one
constructed in the previous proof.
 
The representation $\rho_v$ was constructed so that $\rho_v(t)$
was a minimal lift of $\overline B$. Because $\ell$ does not divide
$\#\rhobar(I_v)$, the matrix $\overline B$ is completely
reducible. So the $V_{i,u}$ in the definition of minimal lift are
$1$-dimensional. Therefore $\rho_v(t)$ is completely reducible
and of finite order prime to~$\ell$. 

The universal property of $(\rho_v',R_v')$ shows that there is a
morphism $R'\to R_v$ which induces $\rho_v$ form $\rho_v'$. Because it
is an isomorphism on mod $\ell$ tangent spaces, the morphism is
surjective. Since both rings are smooth of the same dimension it must
be bijective. This shows that the two deformations agree.}
\end{remark}

\subsection{Proof of Theorem~\ref{OnSerreConj} and the key %
lemma}\label{KeyLemSec}

\begin{Proofof}{Lemma~\ref{keylemma}} We first prove

{\bf Claim 1:} There exists a finite set $T''$ of $R$-places of type
(II) and for each $R$-place $v\in T''$ a choice of eigenvalue
$\lambda_v$, as in the definition, such that 
\begin{equation}\label{FirstOrth} 
H^1_{\{L_v\}}(T'\cup T'',\ad)\cap H^1(\Gal(E(\zeta_\ell)/K),M_n^0(\BF))=0,
\end{equation}
\begin{equation}\label{SecondOrth} 
H^1_{\{L^\perp_v\}}(T'\cup T'',\ad(1))\cap
H^1(\Gal(E(\zeta_\ell)/K),M_n^0(\BF)(1))=0.
\end{equation}
We only give the proof of (\ref{FirstOrth}), the one of
(\ref{SecondOrth}) being analogous. 

Let $\sigma_1,\ldots,\sigma_s$ be the different $R$-classes of type
(II). For $\sigma_i$, let $\lambda_{i,j}$, $j=1,\ldots,m_i$, be the
list of eigenvalues in $\BF$ of multiplicity $2$ (in the characteristic
polynomial). Pick unramified places $v_{i,j}$, $i=1,\ldots,s$,
$j=1,\ldots,m_i$ such that $\Frob_{v_{i,j}}=\sigma_i$ for all $i,j$. 
Let $(\CC_{v_{i,j}},L_{v_{i,j}})$ be the deformation problem defined
in Section~\ref{AuxPrimes} for the pair $(v_{i,j},\lambda_{i,j})$
and let $T'':=\{v_{i,j}:i=1,\ldots,s,j=1,\ldots,m_i\}$. We
consider the following commuting diagram:
$$\xymatrix @C-1.1pc {H^1(\Gal(E(\zeta_\ell)/K),\ad)\cap 
H^1_{\{L_v\}}(T',\ad) \ar[r]\ar[d]&
\prod\limits^{\phantom{v\in T''}}_{v\in T''}\ar[d]
\!\!H^1(\langle\sigma_v\rangle,\ad_{\sigma_v})\\ 
H^1_{\{L_v\}}(T',\ad)\ar[r]&
\prod\limits^{\phantom{v\in T''}}_{v\in T''} \!\!H^1(G_{v},\ad)/L_{v}\rlap{.}
}$$

The kernel of the bottom row is $H^1_{\{L_v\}}(T'\cup T'',\ad)$. 
By assumption there are sufficiently many $R$-classes for
$(\CC_v,L_v)_{v\in S}$, and so the top horizontal arrow is injective.
For each $i$, the image of the top left term in 
$\prod_{j=1,\ldots,m_i}H^1(\langle\sigma_{v_{i,j}}\rangle,\ad_{\sigma_{v_{i,j}}})$
is diagonal. Therefore we may replace the top right term by $\prod_i
H^1(\langle{\sigma_i}\rangle,\ad_{\sigma_i})$ and still retain the
injectivity of the top horizontal map.
Below we show that the induced right vertical arrow 
\begin{equation}\label{NewArrowIsInj}\prod_i
H^1(\langle{\sigma_i}\rangle,\ad_{\sigma_i})\to 
\prod^{\phantom{v\in T''}}_{v\in T''} \!\!H^1(G_{v},\ad)/L_{v} 
\end{equation}
is injective. An easy diagram chase then shows that the intersection
of $H^1_{\{L_v\}}(T'\cup T'',\ad)$ and
$H^1(\Gal(E(\zeta_\ell)/K),\ad)$ is zero, as desired.

The injectivity of (\ref{NewArrowIsInj}) may be verified on the
morphisms
$$ H^1(\langle{\sigma_i}\rangle,\ad_{\sigma_i})\longto
\prod_{j=1,\ldots,m_i} H^1(G_{v_{i,j}},\ad)/L_{v_{i,j}},$$
individually. The representation $\ad_{\sigma_i}$ itself is a direct sum of adjoint representations
on $2\times 2$ blocks of matrices of trace zero. Furthermore, the image of
$H^1(\langle{\sigma_i}\rangle,\ad_{\sigma_i})$ lies in 
$H^1_{\unr}(G_{v_{i,j}},\ad)/L_{v_{i,j},\lambda_{i,j},\unr}$, and
$H^1_{\unr}(G_{v_{i,j}},\ad)$ itself breaks up into a direct sum over
pieces corresponding to the rational canonical form of
$\rhobar(\sigma_i)$. So one is reduced to consider a single $2\times
2$ block, and we may assume $n=2$ and $m_i=1$. But then
$L_{v,\unr}=0$ ($v=v_{i,j}$) and the map 
$$H^1(\langle{\sigma_i}\rangle,\ad)\longto 
H^1_{\unr}(G_{v_{i,j}},\ad)$$
is simply given by inflation and is thus injective. Note finally that
(\ref{FirstOrth}) is preserved under adding further $R$-primes
to~$T''$. We have thus proved Claim~1.

\smallskip

By enlarging $T'$ if necessary, we assume from now on that 
(\ref{FirstOrth}) and (\ref{SecondOrth}) hold with $T''=\emptyset$.
We now induct on the dimension of $ H^1_{\{L_v^\perp\}}(T',\ad)$ and
assume it contains a non-zero cocycle $\phi$. By the formula in
Remark~\ref{WilesFormula} and our assumptions, the space
$H^1_{\{L_v\}}(T',\ad)$ contains a non-zero cocycle $\psi$, as well. 

{\bf Claim 2:} There exists $w\in X\setminus T'$ and an admissible
pair $(\CC_w,L_w)$ compatible with $\eta$ such that the following
hold:
\begin{itemize}
\item[(i)] $n-1=\dim L_w=\dim L_{w,\unr}+1$,
\item[(ii)] $\phi$ does not map to zero in $H^1(G_w,\ad(1))/L_w^\perp$ and
\item[(iii)] the space $(H^1_\unr(G_w,\ad)+L_w)/L_w$ lies in the image of the
morphism $$H^1_{\{L_v\}}(T',\ad)\to H^1(G_w,\ad)/ L_w.$$
\end{itemize}
Suppose for a moment that we have proved the above claim and let
$T'':=T'\cup\{w\}$. The argument given in \cite{taylor}, proof of
Lemma~1.2, then shows, by using (i),(iii), that
$$H^1_{\{L_v\}_{v\in S\cup T'}}(T',\ad)=
H^1_{\{L_v\}_{v\in S\cup T'}\cup \{L_w+H^1_\unr(G_v,\ad)\}}(T'',\ad),$$
and by using (ii), that 
$$H^1_{\{L_v\}_{v\in S\cup T''}}(T'',\ad)\into 
H^1_{\{L_v\}_{v\in S\cup T'}\cup \{L_w+H^1_\unr(G_v,\ad)\}}(T'',\ad).$$
is a proper containment, so that the proof of the lemma is completed.

The above claim is clearly implied by the following {\bf Claim 3},
which we prove below: There exists
$w\in X\setminus T'$ and an admissible pair $(\CC_w,L_w)$ compatible
with $\eta$ such that (i) above holds and furthermore
\begin{itemize}
\item[(ii')] $\phi$ does not map to zero in 
$H^1_\unr(G_w,\ad(1))/L_{w,\unr}^\perp$, and
\item[(iii')] $\psi$ does not map to zero in  
$H^1_\unr(G_w,\ad)/L_{w,\unr}\,(\cong\BF)$.
\end{itemize}

To prove Claim~3, note that 
conditions (\ref{FirstOrth}) and (\ref{SecondOrth}) imply that the
cycles $\psi$ and $\phi$ restrict to non-zero homomorphisms
$$\phi\!:G_{E(\zeta_\ell)}\to (\ad)(1)\quad\mathrm{and}\quad \psi\!:G_{E(\zeta_\ell)}\to \ad.$$
Let $E_\phi$ and $E_\psi$ be the fixed fields of the respective
kernels. Depending on whether the cyclotomic character $\chi$ is
trivial, they may or may not be equal. The induced morphisms on
$\Gal(E_\phi/E(\zeta_\ell))$ and $\Gal(E_\psi/E(\zeta_\ell))$, 
respectively, are equivariant for $\Gal(E(\zeta_\ell)/K)$. Because
$\ad$ is irreducible as a $\BF_\ell[\image(\rhobar)]$-module,
the morphisms $\phi$, $\psi$ are bijective. Thereby the group
$V:=\Gal(E_\psi E_\phi/E(\zeta_\ell))$ may be regarded as an
$\BF[\Gal(E(\zeta_\ell)/K)]$-module which surjects onto $\ad$
and~$\ad(1)$.

Let now $\sigma\in \Gal(E(\zeta_\ell)/K)$ be an $R$-class.
Recall that the subspaces $L_{\unr,\sigma}$ of $\ad/(\sigma-1)\ad$,
and $L^\perp_{\unr,\sigma}$ of $\ad(1)/(\sigma-1)\ad(1)$, defined after
Corollary~\ref{CorOnSigmaAndV}, are of codimension one with respect to
$\BF$. Define 
$\widetilde L_{\unr,\sigma}$ and $\widetilde L^\perp_{\unr,\sigma}$ as the
corresponding $1$-codimensional subspaces in $\ad$ and $\ad(1)$,
respectively.  Each of the conditions 
$$\psi(\tau)\in -\psi(\tilde\sigma)+ \widetilde L_{\unr,\sigma}
\quad\mathrm{and}\quad \phi(\tau)\in  -\phi(\tilde\sigma)+
\widetilde L_{\unr,\sigma}^\perp,$$
$\tau\in V=\Gal(E_\psi E_\phi/E(\zeta_\ell))$, determines a hyperplane in
the $\BF$-vector space $V$. As
we assumed $|\BF|>2$, the join of these two 
hyperplanes cannot span all of $V$, and hence there exists $\tau\in
\Gal(E_\psi E_\phi/E(\zeta_\ell))$ which lies on neither. We fix such a $\tau$ and 
define~$\tilde\tau:=\tau\sigma$.  

The extension $E_\psi E_\phi$ is Galois over $K$, and therefore by
the \v{C}ebotarev density theorem, we can choose a place $w$ in
$X\setminus T'$ such that $\Frob_w=\tilde\tau$. Take $\CC_w$ and $L_w$
as constructed in Section~\ref{AuxPrimes}, so that by
Proposition~\ref{MainPropOnAuxPrimes}, condition (i) is satisfied and $(\CC_w,L_w)$ is
compatible with $\eta$. Condition 
(ii') is satisfied, since the image of $\phi$ in $H^1_\unr(G_w,\ad)$
is given by $\phi(\Frob_w)=\phi(\tilde\tau)=\phi(\tau)
+\phi(\tilde \sigma)\in \ad/(s-1)\ad$, and hence does
not lie in $L_{\sigma,\unr}$. Assertion (iii') is shown in the same
way, and whence the proof of the Lemma~\ref{keylemma} is completed. 
\end{Proofof}

\begin{Proofof}{Theorem~\ref{OnSerreConj}}
By enlarging $X$ if necessary, we may assume that $\rhobar$ ramifies
at all places of $S$. Using Propositions~\ref{FirstPropOnRamV}
and~\ref{SecondPropOnRamV}, there exist locally admissible pairs
$(\CC_v,L_v)_{v\in S}$ compatible with $\eta$ for which one has $\dim
L_v=h^0(G_v,\ad)$ and such that the conductor (at $v$) of any lift of
type $\CC_v$ is the same as that of~$\rhobar_{v}$ . 

We claim that if $\rhobar$ admits sufficiently many $R$-classes, then
it admits sufficiently many $R$-classes for $(\CC_v,L_v)_{v\in S}$. We
only verify the first of the two conditions. For this consider the diagram:
$$\xymatrix{H^1(\Gal(E(\zeta_\ell)/K),\ad)\ar[d]\ar[r] &
\prod^{\phantom{v\in S}}_{v\in S} H^1(\rhobar(I_v),\ad)\ar[d]\\
H^1_{\{L_v\}}(\emptyset,\ad)\ar[r] &\prod^{\phantom{v\in S}}_{v\in S}
H^1(G_{v},\ad)/L_{v}&\rlap{.}\\}$$
To prove the claim, it suffices to show that the right vertical arrow
is injective. By the propositions quoted above, we have
$L_v=H^1_\unr(G_v,\ad)$ for $v\in S$ and thus by the inflation
restriction sequence the morphism $H^1(G_{v},\ad)/L_{v}\into H^1(I_v,\ad)$ is
a monomorphism. Therefore it suffices to show for each $v$ that 
$$ H^1(\rhobar(I_v),\ad)\longto H^1(I_v,\ad)$$
is injective. This is clear since it is an inflation map.

Applying Lemmas~\ref{FirstObs} and~\ref{keylemma}, the assertions of
theorem are now straightforward.
\end{Proofof}

\noindent {\it Addresses of the authors}: 

\noindent GB: Institut f\"ur experimentelle Mathematik, Universit\"at
Duisburg-Essen, Standort Essen, Ellernstrasse 29, 45326 Essen, Germany. \\
e-mail address: {\tt boeckle@exp-math.uni-essen.de}

\noindent CK: Dept of Math, University of Utah, 155 South 1400 East,
Salt Lake City, UT 84112.\\ 
e-mail address: {\tt shekhar@math.utah.edu}
 
\noindent School of Mathematics, 
TIFR, Homi Bhabha Road, Mumbai 400 005, INDIA.\\ 
e-mail addresses: {\tt shekhar@math.tifr.res.in}
 

\begin{thebibliography}{DDT}\itemsep=0pt

\bibitem[B\"o1]{boeckle2} B\"ockle, G., {\it Lifting mod $p$ representations to
characteristic $p^2$}, J.\ Number Theory {\bf 101} (2003), no.\ 2, 310--337.

\bibitem[B\"o2]{boeckle1} B\"ockle, G., {\it Deformations and the rigidity
method}, preprint available under 
{\sf http://www.math.ethz.ch/\~{}boeckle/}.

\bibitem[CPS]{cps} Cline, E., Parshall, B., Scott, L., {\it Cohomology
  of finite groups of Lie type I}, Publ.\ Math.\ IHES {\bf 45} (1975),
  169--191.  

\bibitem[DDT]{ddt} Darmon, H.; Diamond, F.; Taylor, R.\ {\em Fermat's
    last theorem} in `Elliptic curves, modular forms \& Fermat's last
    theorem' (Hong Kong, 1993), 2--140, Internat.\ Press, Cambridge, MA, 1997

\bibitem[dJ]{dejong} de Jong, A. J., {\it Conjecture on arithmetic
    fundamental groups}, Israel J.\ Math.\ {\bf 121} (2001), 61--84. 

\bibitem[Dr]{drinfeld} Drinfeld, V. G., 
{\it Two-dimensional $\ell$-adic representations of the fundamental group of a curve
over a finite field and automorphic forms on ${\rm GL}(2)$}, 
Amer.\ J.\ Math.\ {\bf 105} (1983), 85--114.

\bibitem[Dr1]{drinfeld1}
Drinfeld, V.G., {\it On a conjecture of Kashiwara},  Math. Res. Lett.  8  (2001), 713--728.

\bibitem[Ga]{gaitsgory} D.\ Gaitsgory, {\em On de Jong's conjecture},
  preprint 2004, available under
  http://arXiv.org/abs/math/0402184

\bibitem[Laf]{lafforgue} Lafforgue, L., {\it Chtoucas de Drinfeld et
correspondance de Langlands}, Invent. Math. {\bf 147} (2002), 1--241

\bibitem[NSW]{nsw}
Neukirch, J., Schmidt, A., Wingberg, K., {\it Cohomology of number fields},
Grundlehren 323, Springer-Verlag, Berlin, 2000.

\bibitem[Ra]{ramakrishna} Ramakrishna, R., {\it Deforming Galois
    representations and the conjectures of Serre and Fontaine-Mazur},
    Ann.\ of Math.\ (2)  {\bf 156}  (2002),  no. 1, 115--154. 

\bibitem[Se]{serre} Serre, J.-P., {\it  Sur les
repr\'esentations modulaires de degr\'e 2 de ${\rm
Gal}(\overline{\bf Q}/{\bf Q})$}, Duke Math. J. 54 (1987),
179--230. 

\bibitem[Ta]{taylor} Taylor, R., {\it On icosahedral Artin
representations II}, to appear in American Journal of Math.


\bibitem[Wb]{weibel} Weibel, C., {\it  An introduction to homological
  algebra}, Cambridge Studies in Advanced Mathematics, {\bf 38},
Cambridge University Press, Cambridge, 1994.


\end{thebibliography}
\end{document}